\documentclass[twoside,11pt]{article}%[twoside,12pt]{report}

\usepackage{fancyheadings}%{fancyhdr}
\usepackage{latexsym}
\usepackage{lastpage,engord,graphics}
\usepackage[dvips]{color}
\usepackage{epic}%{pict2e}
\usepackage{eepic}%{pict2e}
\usepackage{rotating}
\usepackage{amssymb}
\usepackage{fancybox}

\makeatletter
\let\savedoddhead=\@oddhead
\def\@oddhead{\hfil {\footnotesize MAPPINGS WITH MAXIMAL RANK}  \hfil\thepage }
\makeatother

\makeatletter
\let\savedevenhead=\@evenhead
\def\@evenhead{\thepage\hfil {\footnotesize C. ABREU-SUZUKI}  \hfil}
\makeatother

\title{Mappings with Maximal Rank}
\author{C. Abreu-Suzuki}
\newcommand{\noin} {\ensuremath{\in\kern-0.77em|}\hspace{0.06in}}
\newcommand{\real} {\ensuremath{I\kern-0.37emR}}
\newcommand{\Na} {\ensuremath{I{\kern-0.37emN}}}
\newcommand{\Z} {\ensuremath{Z\kern-0.40emZ}} % aqui!!!
\newcommand{\rn} {\ensuremath{I\kern-0.37emR^{n}}}
\newcommand{\ra}[1] {\ensuremath{I\kern-0.37emR^{{#1}}}}
\newcommand{\rnz} {\ensuremath{{I\kern-0.37emR}^{n}\setminus\{0\}}}

\newcommand{\pfe} {\mbox{ } \hfill$\Box$}

\newcommand{\pf} {\mbox{\it Proof.\ }}
\newcommand{\cin} {\ensuremath{{\cal C}^{\infty}}}
\newcommand{\ci} {\ensuremath{{\cal C}^{1}}}

\newtheorem{def21}{Definition}[section]
\newtheorem{def22}[def21]{Definition}
\newtheorem{def23}[def21]{Definition}

\newtheorem{lem31}{Lemma}[section]
\newtheorem{lem32}[lem31]{Lemma}
\newtheorem{def313}[lem31]{Definition}
\newtheorem{def33}[lem31]{Definition}
\newtheorem{prop34}[lem31]{Proposition}
\newtheorem{prop35}[lem31]{Proposition}
\newtheorem{lem36}[lem31]{Lemma}

\newtheorem{thm412}{Theorem}[subsection]
\newtheorem{thm413}[thm412]{Theorem}

\newtheorem{thm421}{Theorem}[subsection]
\newtheorem{cex422}[thm421]{Counterexample}
\newtheorem{thm423}[thm421]{Theorem}
\newtheorem{cex424}[thm421]{Counterexample}
\newtheorem{cex425}[thm421]{Counterexample}
\makeindex
\begin{document}
  
   \date{\empty}

%\maketitle
\begin{center}
\textbf{\Large MAPPINGS WITH MAXIMAL RANK}
\end{center}

\begin{center}
{C. ABREU-SUZUKI}
\end{center}

%%%\thispagestyle{fancy}

%\markboth{\protect\footnotesize MAPPINGS WITH MAXIMAL RANK}
%           {\protect\footnotesize C. ABREU-SUZUKI}

%   \noindent{\bf Abstract}
%\begin{abstract}
%Let $\pi: M\rightarrow B$ be an onto maximal rank map or a 
%Riemannian submersion between Riemannian manifolds $M$ and $B$. 
%Initially, we prove necessary and sufficient conditions for any 
%fiber $F$ to be roughly isometric to $M$. Then, we prove 
%necessary and sufficient conditions for $\pi$ to be a rough isometry.  
%As a corollary $M$ is roughly isometric to $F\times B$.
%\end{abstract}

   \pagenumbering{arabic}

\section{Introduction}
%\markboth{\protect\small Mappings with Maximal Rank}
%           {\protect\small Introduction}
%\markright{\protect\small Introduction}
\protect\label{ls1}
\thispagestyle{plain}

We study Maximal Rank Maps and Riemannian Submersions
%\hspace{0.3in} 
$\pi: M\rightarrow B$, where $M$ and $B$ are Riemannian manifolds.

As essential tools in this work we are interested in  equivalence 
relations between non-compact Riemannian manifolds given by 
Rough Isometries, a concept first introduced by 
M. Kanai~\cite{MK1}.

Motivated by O'Neill~\cite{BON} we investigated the question: 
when  does a maximal rank map differ only by a rough isometry of 
$M$  from the simplest type of Riemannian submersions, the 
projection $p_{B}: F\times B \rightarrow B$ of a Riemannian 
product manifold on one of its factors.
Firstly, for Riemannian submersions $\pi: M\rightarrow B$ 
we show that, if the base manifold $B$ is compact and connected, 
then the fibers $F$ can be roughly isometrically immersed into $M$, 
and thus, $M$ is roughly isometric to the product $F\times B$ of any 
fiber and the base space [{\bf Theorem~\ref{lthm412}}]. 
When $B$ is noncompact, connected and complete, and $diam(F)$ is 
uniformly bounded, the Riemannian Submersion $\pi$ is a rough 
isometry, and thus, if a fixed fiber $F$ is 
compact then $M$ is roughly isometric to the product $F\times B$ 
of that fiber and the base space [{\bf Theorem~\ref{lthm413}}].  
Secondly, for onto maximal rank maps that are not necessarily 
submersions, by adding control on the length of 
horizontal vector lifts we have the same consequences 
[{\bf Theorem~\ref{lthm421}}, {\bf Theorem~\ref{lthm423}}]. 
We provide {\it Counterexamples} in section~\ref{ls42}
to show that the assumptions made are necessary conditions.

The paper begins with background.

\section{Rough Isometries and Riemannian Submersions}
%\markboth{\protect\small Mappings with Maximal Rank}
%         {\protect\small Rough Isometries and Riemannian Submersions}
%\markright{\protect\small Rough Isometries and Riemannian Submersions}
\protect\label{ls2}

In this section we define some notation and provide some 
definitions according to M. Kanai~\cite{MK1} %,~\cite{UA} 
and O'Neill~\cite{BON}.

We will be interested in  equivalence relations given by 
rough isometries, a concept first introduced in~\cite{MK1}.

\begin{def21} 
\protect\label{ldef21}
A map $\varphi : M\rightarrow N$, between two metric spaces
$(M,\delta)$ and $(N,d)$, not necessarily continuous, is
called a \underline{rough isometry}, if it satisfies the following two
axioms:
\begin{description}
\item[(RI.1)] There are constants $A\geq 1, C\geq 0$, such that, \\
\[
\frac{1}{A}\delta(p_{1},p_{2}) - C \leq
d(\varphi(p_{1}),\varphi(p_{2}))
\leq A
\delta(p_{1},p_{2}) + C, \hspace{0.2in} \forall p_{1},p_{2}\in M 
\]
\item[(RI.2)] The set $Im\varphi := \{ q = \varphi(p), \forall p \in
M\}$ is
\underline{full} in $N$, i.e.  
\[
\exists\varepsilon >0: N = B_{\varepsilon}(Im\varphi) =
\{q \in N: d(q,Im\varphi)<\varepsilon\}
\] 
In this case we say that $Im\varphi$ is \underline{$\varepsilon$-full} 
in $N$.
\end{description}
\end{def21}

One can easily show that if $\varphi : M\rightarrow N$ and 
$\psi : N\rightarrow M$ are rough isometries, then  the composition 
$\psi\circ\varphi : M\rightarrow M$ is also a rough isometry.
 
We will denote by $\varphi^{-} : N \rightarrow M$ a  \underline{rough
inverse} of $\varphi$, defined as follows: for each $q\in N$, choose 
$p\in M$ so that $d(\varphi(p),q)<\varepsilon$, and define 
$\varphi^{-}(q):= x$. We point out here that 
such a $p$ exists because of the condition {\bf (RI.2)}.  $\varphi^{-}$
is a rough isometry such that both 
$\delta(\varphi^{-}\circ\varphi(p),p)$ and 
$d(\varphi\circ\varphi^{-}(q),q)$
are bounded in $p\in M$ and in $q\in N$, respectively.

We refer to O'Neill~\cite{BON} for the properties of 
Riemannian submersions. We start recalling  their definition.

Let $M^{m}$ and $B^{n}$ be Riemannian manifolds with dimensions 
$m$ and $n$, respectively, where $m\geq n$. 

\begin{def22} 
\protect\label{ldef22}
A map $\pi: M\rightarrow B$ has maximal rank $n$ if the derivative map
$\pi_{\ast}$ is surjective.
\end{def22}

According to~\cite{BON}, a tangent vector on $M$ is said to be
\underline{vertical} if it is tangent to a fiber, \underline{horizontal} 
if it is orthogonal to  a fiber. A vector field on $M$ is 
{vertical} if it is always tangent to  fibers, {horizontal} if it is 
always orthogonal to fibers.

Since the derivative map $\pi_{\ast}{x}$ of $\pi$ is surjective 
for all $x\in M$, its rank is maximal. We can define the 
projections of the tangent space of $M$ onto the subspaces of 
vertical and horizontal vectors, which we will denote respectively 
by $(VT)_{x}$ and $(HT)_{x}$ for each $x\in M$. 
In that case, we can decompose each tangent space to $M$ into a direct
orthogonal sum $T_{x}M=(VT)_{x}\oplus (HT)_{x}$.

\begin{def23} 
\protect\label{ldef23}
A \underline{Riemannian} \underline{submersion} $\pi: M\rightarrow
B$ is an onto mapping satisfying the following two axioms:
\begin{description}
\item[(S.1)] $\pi$ has maximal rank;
\item[(S.2)]$\pi_{\ast}$ preserves lengths of horizontal vectors.
\end{description}
\end{def23}

\section{Long Curves and Their Lifts}
%\markboth{\protect\small Mappings with Maximal Rank}
%           {\protect\small Long Curves and Their Lifts}
%\markright{\protect\small  Long Curves and Their Lifts}
\protect\label{ls3}

Here we begin with background from O'Neill~\cite{BON} and continue 
with an investigation of curves and their lifts.

Let $\pi: M\rightarrow B$ denote an onto mapping with maximal rank $n$
between Riemannian manifolds $M^{m}$ and $B^{n}$ with  $m\geq n$.

From the maximality of the rank of the onto mapping $\pi$ we have the
unique horizontal vector property:

\begin{lem31}
\protect\label{llem31}
 Let $b\in B$. Given any $w\in T_{b}B$ and $x\in M$
satisfying $\pi(x)=b$, there exists  a unique horizontal vector
 $v\in T_{x}M$ which is $\pi$-related to $w$, i.e. satisfying  $v\in
(HT)_{x}$ and  $(\pi_{\ast})_{x}(v)=w$.%doesn't use S.2
\end{lem31}

If, in addition, one has control from below over the length  
of horizontal vectors, then one has contro; from below over the 
distance in $M$. This is the essence of  the following Lemma.

\begin{lem32}
\protect\label{llem32}
Assume that $M$ and $B$ are both connected and geodesically complete.
Let $x,x'\in M$, $\Gamma_{min}\subset M$ be a minimal geodesic joining
$x$ to $x'$, and $\gamma_{min}\subset B$ be a minimal geodesic joining
$\pi(x)$ to $\pi(x')$.
Suppose that for all $b \in B$ and for all
$x\in F_{b}$ there exist constants $\alpha\geq 1$ and $\beta >0$,
both independent of $b$ and $x$, such that
\begin{equation}
\protect\label{lequ47}
\frac{1}{\alpha} ||w||_{B} - \beta \leq ||v||_{M}
\end{equation}
for all $w\in T_{b}B$,  where $v$ is the unique horizontal lift of
$w$ through $x$ that we assume satisfies $||v||_{M}\leq 1$, and
$||\hspace{0.15in}||_{M}$, $||\hspace{0.15in}||_{B}$ denote the
inner product on $TM$ and $TB$, respectively.

Then, \hspace{0.1in}
$
d_{M}(x,x')=
\ell(\Gamma_{min})\geq \frac{1}{\alpha} \ell(\gamma_{min})-\beta
=\frac{1}{\alpha}d_{B}(\pi(x),\pi(x'))-\beta
$
\end{lem32}
\pf\hspace{0.1in}
Without loss of generality, we may assume that the horizontal lift 
$v\in (HT)_{x}$ of $w$ satisfies $||v||_{M}\leq 1$, and that  
both parametrizations of $\Gamma_{min}$ and $\gamma_{min}$ are 
defined in the interval $[0,1]$.

We may write 
\[
\Gamma'_{min}(t)=\Gamma'_{V}(t)\oplus\Gamma'_{H}(t)\in
T_{x}M=(VT)_{x}\oplus (HT)_{x}, \forall t\in [0,1]
\]

Notice that by {\bf Lemma~\ref{llem31}},
$\Gamma'_{H}(t)\in T_{\Gamma_{min}(t)}M$ is the unique horizontal
vector which is $\pi-$related to
$\frac{d}{dt}\left(\pi\circ\Gamma_{min}\right)(t)$,  for each
$t\in [0,1]$. Assume that $\Gamma$ is parametrized proportionally 
to arclength, and that $||\Gamma'_{H}(t)||_{M}\leq 1$.

We have,
\begin{eqnarray*}
\lefteqn{d_{M}(x,x') =
\ell(\Gamma_{min})=\int_{0}^{1}||\Gamma'_{V}(t)\oplus
\Gamma'_{H}(t)||_{M} dt \geq \int_{0}^{1}||\Gamma'_{H}(t)||_{M} dt
\stackrel{(\ref{lequ47})}{\geq}} \\
& \geq & \frac{1}{\alpha} \int_{0}^{1}\left|\left|
\frac{d}{dt}\left(\pi\circ\Gamma_{min}\right)(t)\right|\right|_{M} dt
-\beta = \frac{1}{\alpha}\ell(\pi\circ\Gamma_{min}) - \beta
\stackrel{\mbox{\scriptsize min.geod.}}{\geq}\\
& \geq & \frac{1}{\alpha} \ell(\gamma_{min})-\beta
\stackrel{\mbox{\scriptsize dist.}}{=}
\frac{1}{\alpha}d_{B}(\pi(x),\pi(x'))-\beta
\end{eqnarray*}
which concludes the Lemma.

\pfe

In what follows lifts of curves are defined.

\begin{def313}
\protect\label{ldef313}
Let $\gamma:[t_{1},t_{2}]\rightarrow B$ be a smooth embedded curve 
in $B$ and $\Gamma:[t_{1},t_{2}]\rightarrow M$ be any curve in $M$ 
satisfying $\pi\circ\Gamma=\gamma$.  The curve $\Gamma$ is called a 
\textbf{lift} of $\gamma$.

If in addition, $\Gamma$ is horizontal, 
i.e., ${\Gamma}^{\prime}(t)\in (HT)_{\Gamma(t)}, 
\forall t\in [t_{1},t_{2}]$, where 
$\Gamma(t_{1})=x_{0}\in M$ with $\gamma(t_{1})=\pi(x_{0})$, the curve 
$\Gamma$ is called a \textbf{horizontal lift} of $\gamma$ 
through $x_{0}$.
Recall that the horizontal lift of a curve in $B$, through a 
point $x_{0}\in M$ is unique.
\end{def313}

\vspace{0.1in}

Next, we define long curves.

\begin{def33}
\protect\label{ldef33}
Let $\beta$ be any positive constant. A smooth embedded curve 
$\gamma:[t_{1},t_{2}]\rightarrow B$ is said to be a 
\textit{$\beta$-long curve} if 
$\inf_{t_{1}\leq t\leq t_{2}}||\gamma'(t)||\geq\beta$. In that case, 
$\ell(\gamma)\geq \displaystyle{\int}_{t_{1}}^{t_{2}}||\gamma'(t)|| dt 
\geq \beta (t_{1}-t_{2})$. We say that a curve $\gamma$ is simply a 
long curve if it is a $\beta$-long curve for some constant $\beta >0$.
\end{def33}

Let $\gamma:[t_{1},t_{2}]\rightarrow B$ denote a smooth embedded curve
and let $\Gamma:[t_{1},t_{2}]\rightarrow M$ denote a lift of $\gamma$.

In the next two  Propositions, under control from above (below)  on the
derivative of the maximal rank mapping $\pi$,  we have  control  from 
below (above) over the length of any lift of a curve. 

For instance,  in
{\bf Proposition~\ref{lprop34}}
for a long curve $\gamma$ in $B$ any of its lift
$\Gamma$ in $M$ cannot be short, and  {\bf Proposition~\ref{lprop35}}
the length of a lift $\Gamma$ of a long curve $\gamma$ is
bounded above by the  length of $\gamma$.

 We denote by $||\hspace{0.1in}||_{M}$ and $||\hspace{0.1in}||_{B}$ the
Riemannian norms in $TM$ and $TB$, respectively.

\begin{prop34}
\protect\label{lprop34}
 Assume there are constants $\alpha\geq 1$ and $\beta >0$ such that,
\begin{equation}
\protect\label{lequ42}
\left||(\pi_{\ast})_{x}v|\right|_{B} \leq
\alpha\left||v|\right|_{M}+\beta
\end{equation}
for all $x\in M$, for all $v \in T_{x}M$ satisfying
$\left||v|\right|_{M}\leq 1$.

If $\gamma$ is any smooth $\beta$-long curve in $B$, then,
\[
\ell(\Gamma)\geq\frac{1}{\alpha}\left[\ell(
\gamma)-\beta (t_{2}-t_{1})\right]>0
\]
where $\ell(\Gamma)$ and $\ell(\gamma)$ denote the lengths of
the curves $\Gamma$ and $\gamma$, respectively.
\end{prop34}
\pf\hspace{0.1in}
First, we choose a parametrization proportional to arc length of
$\Gamma: t\in[t_{1},t_{2}] \rightarrow \Gamma(t)\in M$,  an arbitrary
lift of $\gamma\subset B$.  We may assume without loss of 
generality that $\left||\Gamma'(t)|\right|_{M} \leq 1$.

If we use $v=\Gamma'(t)$ in~(\ref{lequ42}) and 
$\pi\circ\Gamma = \gamma$, we obtain
\begin{eqnarray}
\protect\label{leqn43}
\left||\gamma'(t)|\right|_{B} = 
\left||(\pi_{\ast})_{\Gamma(t)}\Gamma'(t)|\right|_{B} \leq
\alpha\left||\Gamma'(t)|\right|_{M}+\beta, & &
\forall t \in [t_{1},t_{2}]
\end{eqnarray}

Finally, if we integrate~(\ref{leqn43}), we get
\[
\begin{array}{rcl}
\ell(\gamma) & = &
{\displaystyle\int_{t_{1}}^{t_{2}}}\left||\gamma'(t)|\right|_{B} dt
\leq  \alpha{\displaystyle\int_{t_{1}}^{t_{2}}}
\left||\Gamma'(t)|\right|_{M} dt +
\beta{\displaystyle\int_{t_{1}}^{t_{2}}} dt = \\
& & \\
& = & \alpha\cdot\ell(\Gamma) + \beta(t_{2}-t_{1})\Rightarrow\\
& & \\
&  \Rightarrow & \ell(\Gamma)\geq
\frac{1}{\alpha}\left[\ell( \gamma)-\beta(t_{2}-t_{1})\right]>0
\end{array}
\]
which proves the proposition.  

That the second hand side of the last inequality above is  
positive 
follows from the assumption that $\gamma$ is a long curve. 
Therefore, as it can be interpreted from the inequality shown, for a
long curve $\gamma$ any of its  lift $\Gamma$ cannot be short.

\pfe

\begin{prop35}
\protect\label{lprop35} Let $\gamma:[t_{1},t_{2}]\rightarrow B$ denote
a smooth embedded curve and let $\Gamma:[t_{1},t_{2}]\rightarrow M$
denote a lift of $\gamma$.
For \underline{horizontal}$^{(\dagger)}$
vectors $v\in TM$ only,
assume that there is a universal constant $\alpha\geq 1$
such that,
\begin{equation}
\protect\label{lequ46}
\left||(\pi_{\ast})_{x}v|\right|_{B} \geq
\frac{1}{\alpha}\left||v|\right|_{M}-\beta
\end{equation}
for all $x\in M$, for all $v \in T_{x}M\setminus (VT)_{x}
=(HT)_{x}=\left[\ker (\pi_{\ast})_{x}\right]^{\perp}$.

If $\gamma$ is a $\beta$-long curve then,
\[
\ell(\Gamma)\leq \alpha\left[\ell(\gamma)+\beta (t_{2}-t_{1})\right]
\]
where $\ell(\Gamma)$ and $\ell(\gamma)$ denote the lengths of
the curves $\Gamma$ and $\gamma$, respectively.
\end{prop35}

Now, in the proof of {\bf Proposition~\ref{lprop35}} we will need
the following Lemma: for a long curve in $B$, any of its  lift  in $M$
 is non-vertical.

\begin{lem36}
\protect\label{llem36}
Let $\gamma:[t_{1},t_{2}]\rightarrow B$ be a smooth embedded curve in
$B$ and let $\Gamma:[t_{1},t_{2}]\rightarrow M$ be a lift of $\gamma$.
If $\gamma$ is a long curve,
then, $\Gamma$ is non-vertical, i.e. there exists an interval
$[t_{1},t_{2}]$, such that,
\begin{eqnarray*}
\left(\Gamma'(t)\right)_{H}\neq 0, & & \forall t\in [t_{1},t_{2}]
\end{eqnarray*}
where $\Gamma'(t)=\left(\Gamma'(t)\right)_{V}\oplus \left(\Gamma'(t)
\right)_{H}\in T_{\Gamma(t)}M= (VT)_{\Gamma(t)}\oplus (HT)_{\Gamma(t)}$.
\end{lem36}
\pf\hspace{0.1in}
Since $\gamma$ is a smooth, embedded long curve, there exists an interval let us say
$[t_{1},t_{2}]$, for which,
\[
\gamma'(t) \neq 0, \hspace{0.1in} \forall t\in [t_{1},t_{2}]
\]

Moreover, since for all $t\in [t_{1},t_{2}]$ the restriction of the 
derivative map 
\[
(\pi_{\ast})_{\Gamma(t)}\mid_{(HT)_{\Gamma(t)}} 
\mbox{ is an isomorphism, }
\]
we thus obtain
\[
(\pi_{\ast})_{\Gamma(t)} \left(\Gamma'(t)\right)_{H} =
(\pi_{\ast})_{\Gamma(t)}
    \left\{\left(\Gamma'(t)\right)_{V}\oplus \left(\Gamma'(t)\right)_{H}
    \right\} =
(\pi_{\ast})_{\Gamma(t)} \left\{\Gamma'(t)\right\}=
\gamma'(t) \neq 0
\]
\[
\stackrel{isom.}{\Longrightarrow}
\left(\Gamma'(t)\right)_{H}\neq 0
\]
for all $t\in [t_{1},t_{2}]$, and thus $\Gamma$ is non-vertical.

\pfe

\pf\textit{of Proposition~\ref{lprop35}}\hspace{0.1in}  We first 
notice that because $\gamma$ is a long
curve, by {\bf Lemma~\ref{llem36}}, $\Gamma$ is non-vertical.

If we use the horizontal vector $v=\Gamma'(t)$ in~(\ref{lequ46}), 
we may write
\begin{eqnarray*}
0\stackrel{Prop~\ref{lprop35}(\dagger)}{\neq}
\left||(\pi_{\ast})_{\Gamma(t)}\Gamma'(t)|\right|_{B} \geq
{\displaystyle\frac{1}{\alpha}}\left||\Gamma'(t)|\right|_{M}-\beta
, & & \forall t \in [t_{1},t_{2}]
\end{eqnarray*}
and using $\pi\circ\Gamma=\gamma$ in the above inequality, we obtain
\begin{eqnarray*}
0\neq\left||\gamma'(t)|\right|_{B} \geq
{\displaystyle\frac{1}{\alpha}}\left||\Gamma'(t)|\right|_{M}-\beta
, & & \forall t \in [t_{1},t_{2}]
\end{eqnarray*}
which in turn implies that
\begin{eqnarray*}
\left||\Gamma'(t)|\right|_{M} \leq
\alpha\left(\left||\gamma'(t)|\right|_{B}+\beta\right)
, & & \forall t \in [t_{1},t_{2}]
\end{eqnarray*}

Finally, integrating the above inequality gives us
\[
\begin{array}{rcl}
\ell(\Gamma) &=& {\displaystyle\int_{t_{1}}^{t_{2}}}\left||\Gamma'(t)
\right||_{M} dt \leq 
\alpha {\displaystyle\int_{t_{1}}^{t_{2}}}
\left||\gamma'(t)\right||_{B} dt
+ \beta\alpha{\displaystyle\int_{t_{1}}^{t_{2}}} dt =\\
& = & \alpha\cdot\ell(\gamma)+ \beta\alpha(t_{2}-t_{1})=
\alpha\left[\ell(\gamma)+\beta(t_{2}-t_{1})\right]>0
\end{array}
\]
which proves the proposition.

Therefore, as it can be interpreted from the above inequality, the
length of a lift $\Gamma$ of a long curve $\gamma$ is
controlled by above by the  length of $\gamma$.

\pfe

\section{Riemannian Submersions, Maximal Rank Maps and Counterexamples}
%\markboth{\protect\small Mappings with Maximal Rank}
%         {\protect\small Riemannian Submersions, Maximal Rank Maps and 
%Counterexamples}
%\markright{\protect{\small Riemannian Submersions, Maximal Rank Maps and 
%Counterexamples}}
\protect\label{ls4}

In this section we will explore Riemannian submersions and maximal rank 
maps $\pi: M\rightarrow B$ between Riemannian manifolds $M$ and $B$.

Motivated by O'Neill~\cite{BON}, we will investigate this question:  
when does a maximal rank map $\pi: M\rightarrow B$ differ only by a 
rough isometry of $M$ from the simplest type of Riemannian submersions, 
the projection $p_{B}: F\times B \rightarrow B$ of a Riemannian product 
manifold on one of its factors.

\subsection{Riemannian Submersions}
%\markboth{\protect\small Riemannian Submersions, Maximal Rank 
%Maps and Counterexamples}
%           {\protect\small Riemannian Submersions}
%\markright{\protect\small Riemannian Submersions}
\protect\label{ls41}

We first show that, 
if the base manifold $B$ is compact and connected, 
then the fibers $F$ can be roughly isometrically immersed into $M$, 
and thus, $M$ is roughly isometric to the product $F\times B$ of any 
fiber $F$ and the base space $B$ [{\bf Theorem~\ref{lthm412}}]. 
Secondly, when $B$ is noncompact, connected and complete, and 
$diam (F)$ is uniformly bounded, we show that the Riemannian 
submersion $\pi: M\rightarrow B$ is a rough isometry, and thus, if a 
fixed fiber $F$ is compact then $M$ is roughly isometric to the 
product $F\times B$ of that fiber $F$ and the base space $B$ 
[{\bf Theorem~\ref{lthm413}}].

\begin{thm412}
\protect\label{lthm412}
 Let $\pi: M\rightarrow B$ be a Riemannian submersion.
Suppose $B$ is compact and connected, and for each $b\in B$ the
fiber
$\pi^{-1}(b)$ has the induced metric from $(M,d)$. 
Then for each $b\in B$, the inclusion
$\iota:\pi^{-1}(b) \hookrightarrow M$ is a rough isometry. 

In particular, since $B$ is compact, $M$ is roughly isometric to the 
product $\pi^{-1}(b)\times B$. 
\end{thm412}

\begin{thm413}
\protect\label{lthm413}
Let $\pi: M\rightarrow B$ be a Riemannian submersion,
where $B$ is  connected and complete. Suppose that,
for some constant $m>0$, all fibers satisfy the
universal diameter property:
\begin{description}
\item[(UDF)]
$\mbox{diam } \left(\pi^{-1}(b)\right) \leq m < \infty, \forall b\in B$.
\end{description}
Then, $\pi: M\rightarrow B$ is a rough isometry. 

In particular, if for some $b_{0}$ the fiber $\pi^{-1}(b_{0})$ is 
compact, then $M$ is roughly isometric to the product 
$\pi^{-1}(b_{0})\times B$. 
\end{thm413}

Note that 
{\bf Theorem~\ref{lthm412}} is a Corollary of 
{\bf Theorem~\ref{lthm421}}, 
and {\bf Theorem~\ref{lthm413}} is a Corollary of 
{\bf Theorem~\ref{lthm423}}, proven in the next section.

\subsection{Non-Submersions  \hspace{0.0005in} Surjective 
\hspace{0.0005in} Maximal \hspace{0.0005in} Rank \hspace{0.0005in} 
Maps}
%\markboth{\protect\small Riemannian Submersions, Maximal Rank 
%Maps and Counterexamples}
%%         {\protect\small Non-Submersions Surjective Maximal 
%Rank Maps}
%\markright{\protect Non-Submersions Surjective Maximal Rank Maps}
\protect\label{ls42}

In this section, we prove that for onto smooth mappings with maximal 
rank $\pi:M\rightarrow B$,
that are not necessarily submersions, the same results as in 
{\bf Theorem~\ref{lthm412}} and {\bf Theorem~\ref{lthm413}} hold,
as long as we make extra assumptions on the subspaces of horizontal
vectors by adding control from above on the length of 
horizontal vector lifts   
[{\bf Theorem~\ref{lthm421}}, {\bf Theorem~\ref{lthm423}}]. 
{\bf Counterexamples} are provided in this section%~\ref{ls42}
 to show that if any of the assumptions are removed those results 
cease to follow 
[{\bf Counterexample~\ref{lcex422}}, {\bf Counterexample~\ref{lcex424}},
 {\bf Counterexample~\ref{lcex425}}].

\begin{thm421}
\protect\label{lthm421}
Assume that $B$ is compact, and  for each $b\in
B$, the fiber $\pi^{-1}(b) = F_{b}$ is endowed with the  induced
metric from $(M,d)$. Suppose that for all $b\in B$, 
there are constants $\alpha\geq 1$ and $\beta >0$, independent of $b$
such that,  the following inequality holds:
\begin{equation}
\protect\label{leqn416}
||v||_{M} \leq \alpha ||w||_{B} + \beta
\end{equation}
for all $x\in F_{b}$ and $w\in T_{b}B$, where
$v\in (HT)_{x} \subset T_{x}M$ is the horizontal lift of $w$
through $x$.

Then,  for each $b\in B$,  the inclusion  map $\iota:F_{b}
\hookrightarrow M$ is a rough isometry.

In particular, since $B$ is compact, $M$ is roughly isometric 
to the product $\pi^{-1}(b)\times B$. 
\end{thm421}
\pf\hspace{0.1in} We must verify axioms {\bf (RI.1)} and {\bf (RI.2)}
for $\iota$, given any $b\in B$.

Clearly axiom {\bf (RI.1)} holds since each fiber 
has the induced metric.

Let us denote by $(M,d_{M})$ and $(B,d_{B})$ the Riemannian metric
spaces, and let $b\in B$ be fixed.

To verify axiom {\bf (RI.2)}, we need to prove that $M$ is an
$\epsilon$-neighborhood of $\iota(F_{b})\subseteq M$, for some
$\epsilon >0$, i.e. we must find a constant $\epsilon >0$ for which
\[
d_{M}(y,\iota(F_{b}))<\epsilon, \hspace{0.1in} \forall  y \in M
\]

Without loss of generality, we may assume that $B$ is connected, 
otherwise, we can repeat, on each connected component of $B$, the 
argument that will follow.

Now, since $B$ is compact and connected it is also complete.

Thus, for any $y\in M$ there exists a minimal geodesic $\gamma$
joining $\pi(y)$ to $b$, which we will parametrize by 
\[
\gamma: [0,1]  \rightarrow  B,
\gamma(0)  =  \pi(y), \gamma(1) =  b
\]

Since $\gamma$ has a unique horizontal
lift $\Gamma_{y}:[0,1]\rightarrow M$, through $y$, and so $\Gamma_{y}$ 
connects $y$ to the fiber $F_{b}$, we can write,
\begin{eqnarray}
\protect\label{leqn417}
d_{M}(y,\iota(F_{b})) & \stackrel{dist.}{\leq} & \ell(\Gamma_{y}) :=
\int_{0}^{1}||\Gamma'_{y}||_{M} dt
\stackrel{(\ref{leqn416})}{\leq}  \nonumber \\
& \leq & \alpha \int_{0}^{1}
||(\pi_{\ast})_{\Gamma_{y}(t)}\Gamma'_{y}(t)||_{B} dt
+ \beta
=  \nonumber \\
& \stackrel{\pi\circ\Gamma_{y}=\gamma}{=} & \alpha \int_{0}^{1}
||\gamma'(t)||_{B} dt + \beta
= \alpha\cdot \ell(\gamma) + \beta
\end{eqnarray}

Now, by the compactness of $B$,
\[
\mbox{diam } B :=\sup_{b_{1},b_{2}\in B}\{d_{B}(b_{1},b_{2})\}<\infty
\]

Moreover, since $\gamma$ is a minimal geodesic joining $\pi(y)$ to
$b$,
\begin{equation}
\protect\label{lequ418}
\ell(\gamma)=d_{B}(\pi(y),b) \leq \mbox{diam } B
<\infty
\end{equation}

By substituting (\ref{lequ418}) in (\ref{leqn417}), we obtain,
\begin{eqnarray}
\protect\label{leqn419}
d_{M}(y,\iota(F_{b}))\stackrel{(\ref{leqn417})}{\leq} \alpha\cdot
\ell(\gamma) + \beta
\stackrel{(\ref{lequ418})}{\leq} \alpha\cdot
(\mbox{diam } B) + \beta
\end{eqnarray}

Define $\epsilon:= \alpha\cdot (\mbox{diam } B) + \beta$,
which is a positive constant independent of $y$, and also of $b \in B$.

For that choice of $\epsilon$, since $y\in M$ is arbitrary, we see that
(\ref{leqn419}) is exactly axiom {\bf (RI.2)} for the inclusion map 
$\iota:F_{b} \hookrightarrow M$.

\pfe

In what follows, we provide a Counterexample to illustrate how
assumption (\ref{leqn416}) is essential in
{\bf Theorem~\ref{lthm421}}. We show that if (\ref{leqn416}) doesn't
hold for some $b\in B$, then the inclusion map
$\iota:F_{b} \hookrightarrow M$ ceases to be a rough isometry.

\begin{cex422}
\protect\label{lcex422}
We will exhibit $M,B,\pi$ satisfying all the conditions in
\textsf{Theorem~\ref{lthm421}} with the exception of
\textsf{(\ref{leqn416})},
i.e.,

For any given constants $\alpha\geq 1$ and $\beta >0$ there exist
$b\in B$, $w\in T_{b}B$ and $x\in F_{b}$ satisfying:
\begin{eqnarray}
\protect\label{leqn440}
||v||_{M} > \alpha ||w||_{B} + \beta %\ce{1}
\end{eqnarray}
where $v$ is the  horizontal lift of $w$ in $(HT)_{x}\subset T_{x}M$.

In this case, the inclusion map $\iota:F_{b} \hookrightarrow M$ is not
a rough isometry.
\end{cex422}

      \begin{figure}[here]

         %\begin{picture}(100,350) (0,0)
         \begin{picture}(385,330)(0,0)

%\dottedline{2}(0,350)(385,350)
%\dottedline{2}(0,0)(385,0)

\put(5,315){\shadowbox{\Huge M}}

%parallels=horizontal curves

%\qbezier(5,320)(60,340)(115,320)
%\qbezier(5,320)(60,300)(115,320)
\put(60,300){\ellipse{110}{20}}

%\qbezier[90](23,225)(60,240)(97,225)
%\qbezier(23,225)(60,210)(97,225)
\put(60,205){\ellipse{75}{15}}

%\qbezier[90](5,130)(60,150)(115,130)
%\qbezier(5,130)(60,110)(115,130)
\put(60,110){\ellipse{110}{20}}

%meridians=vertical curves
\qbezier(5,110)(40,205)(5,300)

\thicklines
\qbezier(85,101)(55,200)(85,291)
\thinlines

\qbezier(115,110)(80,205)(115,300)

%text
\put(44,271){\ovalbox{\Large $F_{b}$}$\leadsto$}

\put(70,160){$\bullet$}
\put(77,160){$\xi_{s}$}

\put(115,225){{\scriptsize $(x_{1},x_{2},x_{3}) \mapsto
\left(\displaystyle{\frac{x_{1}}{\sqrt{x_{1}^{2}+x_{2}^{2}}}},
\displaystyle{\frac{x_{2}}{\sqrt{x_{1}^{2}+x_{2}^{2}}}}, x_{3}\right)$}}
\put(105,205){\vector(1,0){180}}
\put(155,190){(diffeomorphism)}

\put(60,95){\vector(0,-1){50}}
\put(65,70){\Large $\pi$}

\put(5,40){\shadowbox{\Huge B}}

%\qbezier(18,30)(55,45)(92,30)
%\qbezier(18,30)(55,15)(92,30)

%\put(60,50){\oval(65,32)}
\put(60,30){\ellipse{75}{20}}
\put(73,18){$\bullet$}%\put(73,41){\circle*{5}}
\put(73,8){$b$}

%cylinder
%parallels=horizontal curves

%\qbezier(294,320)(331,340)(368,320)
%\qbezier(294,320)(331,300)(368,320)
\put(331,300){\ellipse{74}{20}}

%\qbezier[90](294,225)(331,240)(368,225)
%\qbezier(294,225)(331,210)(368,225)
\put(331,205){\ellipse{74}{20}}

%\qbezier[90](294,130)(331,150)(368,130)
%\qbezier(294,130)(331,110)(368,130)
\put(331,110){\ellipse{74}{20}}

%meridians=vertical curves
\put(294,110){\line(0,1){190}}
\put(368,110){\line(0,1){190}}
%axis
\thicklines
\dashline[+70]{3}(331,206)(331,304)
\thinlines
\put(331,304){\vector(0,1){30}}%\put(331,230){\vector(0,1){120}}
    \put(335,323){$e_{3}$}
\thicklines
\dashline[+70]{3}(331,206)(359,198)
\put(359,198){\vector(3,-1){26}}%\put(331,226){\vector(3,-1){45}}
    \put(371,182){$e_{2}$}
\thicklines
\dashline[+70]{3}(331,206)(303,198)
\put(303,198){\vector(-3,-1){26}}%\put(331,226){\vector(-3,-1){45}}
    \put(282,182){$e_{1}$}
\thinlines
%text
\put(155,40){$(\tilde{a},\tilde{b},\tilde{c})\mapsto
             (\tilde{a},\tilde{b},0)$}
\put(331,95){\line(0,-1){45}}%{65}}
\put(306,55){\oval(50,50)[br]}
\put(305,30){\vector(-1,0){200}}%\put(331,50){\vector(-1,0){220}}
\put(174,17){(projection)}

%\put(5,4){\line(0,-1){2}}
%\put(6,2){\oval(2,2)[bl]}
%\put(6,1){\vector(1,0){6}}
   \end{picture}

         \caption{The map $\pi: M \rightarrow B$ in
         {\bf Counterexample~\ref{lcex422}}.}
         \label{fig44cex422}
         \index{pictures!Counterexample\ref{lcex422}}
      \end{figure}
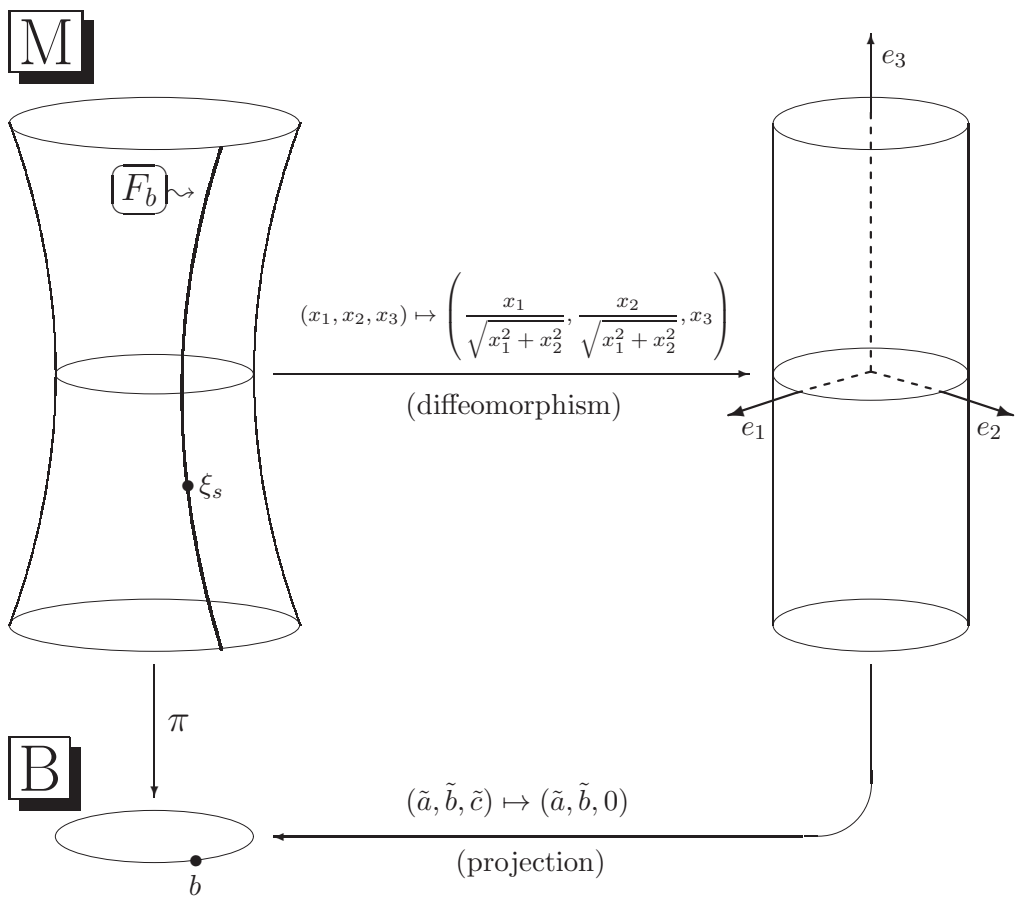

%(see Fig.~{\ref{fig44cex422}})

Let $M$ and $B$ be the following Riemannian manifolds,
\[
M=\{(x_{1},x_{2},x_{3})\in\real^{3}: x_{1}^{2}+x_{2}^{2}=x_{3}^{2}+ 1\}
\]
and the compact unit circle,
\[
B={\cal S}^{1}=\{(x_{1},x_{2},0)\in \real^{3}:  x_{1}^{2}+x_{2}^{2}
= 1\}
\]
where the metrics on $M$ and $B$ are induced by the Euclidean metric
on $\real^{3}$.

Let $\pi: M \rightarrow B$ be defined by,
\[
\pi(x_{1},x_{2},x_{3})=\left(\displaystyle{\frac{x_{1}}{\sqrt{x_{1}^{2}
+x_{2}^{2}}}},
\displaystyle{\frac{x_{2}}{\sqrt{x_{1}^{2}+x_{2}^{2}}}}, 0\right)
\]

Clearly $\pi: M\rightarrow B$ is an
\underline{onto smooth maximal rank} map. 

Firstly, we remark that (\ref{leqn440}) can be verified with a 
series of calculations (c.f.~\cite{CAS}).

Lastly, 
we show that for each $b=(b_{1},b_{2},0)\in B$
the inclusion
$\iota: F_{b}\rightarrow M$ is \underline{not} a rough isometry.

In that direction, we claim that {\bf (RI.2)} fails, i.e.
\[
\forall \epsilon >0, \exists y_{\epsilon}=y_{\epsilon}(\epsilon,b)
\in M, \mbox{ satisfying } d_{M}(y_{\epsilon},\iota(F_{b}))
\geq\epsilon
\]

Let $\gamma$ be a compact connected smooth curve in
$B={\cal S}^{1}$, parametrized by,
\[
\gamma(t) = (\cos(t), \sin(t),0)\in B, \hspace{0.1in} \forall t\in[0,1]
\]
with $b=\gamma(t_{b})$ for some $t_{b}\in[0,1]$.

A generic element in the fiber $F_{b}\subseteq M$ can be described as,
\[
\xi_{r}:=\left(\cos(t_{b})\cdot\sqrt{r^{2}+1},
\sin(t_{b})\cdot\sqrt{r^{2}+1}, r\right) \in M
\]
where  $s\in\real$ is constant.

Thus, the fiber $F_{b}$, where
$b=\gamma(t_{b})=(\cos t_{b},\sin t_{b},0)$, can be described as,
\[
F_{b}=\{\xi_{r}:=(\sqrt{r^{2}+1}\cdot\cos t_{b},
         \sqrt{r^{2}+1}\cdot\sin t_{b},r), r\in \real\}
\]

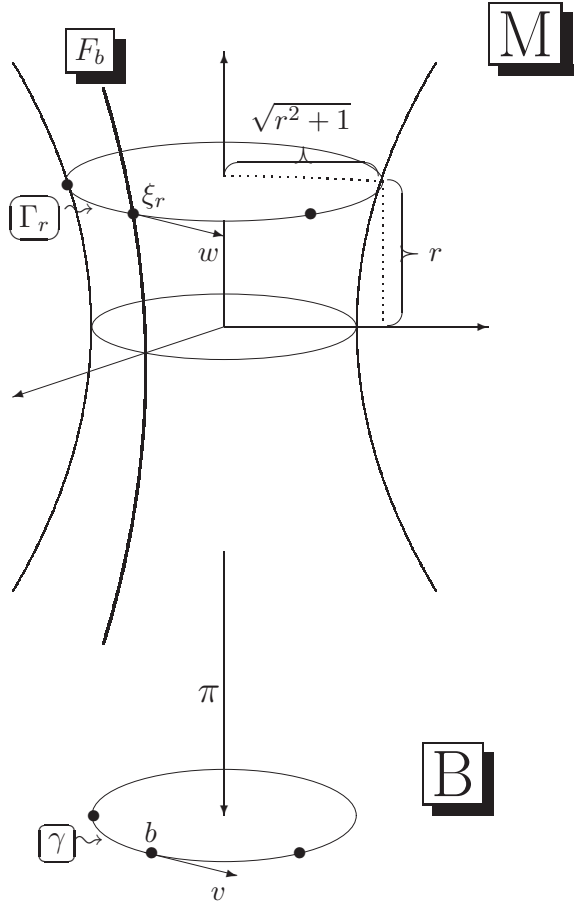
\begin{figure}[here]

         %\begin{picture}(100,350) (-60,10)
            %\includegraphics{piccex422-5.ps}
        \begin{picture}(300,340)(0,0)

%\dottedline{2}(0,350)(300,350)
%\dottedline{2}(0,0)(300,0)

\put(270,320){\shadowbox{\Huge M}}
%parallels=horizontal curves

%\qbezier[90](110,300)(170,320)(230,300)
%\qbezier(110,300)(170,280)(230,300)
\put(170,280){\ellipse{120}{30}}

%\qbezier[90](120,245)(170,260)(220,245)
%\qbezier(120,245)(170,230)(220,245)
\put(170,225){\ellipse{100}{25}}

%meridians=vertical curves
\qbezier(90,125)(149,225)(90,325)

\thicklines
\qbezier(124,105)(156,210)(124,315)
\put(110,317){\shadowbox{\large $F_{b}$}}

\thinlines
\qbezier(250,125)(190,225)(250,325)
%axis
\put(170,225){\line(0,1){40}}
\put(170,282){\vector(0,1){47}}
\put(170,225){\vector(1,0){100}}
\put(170,225){\vector(-3,-1){80}}
%text
\put(108,276){$\bullet$}%\put(110,300){\circle*{5}}
\put(133,265){$\bullet$}%\put(135,288){\circle*{5}}
\put(200,265){$\bullet$}%\put(200,288){\circle*{5}}

\put(135,268){\vector(4,-1){35}}
\put(160,250){$w$}
\put(139,272){$\xi_{r}$}
\put(89,263){\ovalbox{\large $\Gamma_{r}$}}
\put(110,267){$\leadsto$}

\put(180,300){$\sqrt{r^{2}+1}$}
\put(196,287){$\curlywedge$}
\qbezier[20](170,282)(200,282)(230,280)
\put(199,282){\oval(58,10)[t]}

\put(236,250){$\succ r$}
\qbezier[20](230,280)(230,252)(230,225)
%\put(220,271){\vector(0,1){30}}
\put(232,253){\oval(10,55)[r]}

\put(245,40){\shadowbox{\Huge B}}

\put(170,140){\vector(0,-1){100}}
\put(160,83){\Large $\pi$}

%\qbezier(38,30)(75,45)(112,30)
%\qbezier(38,30)(75,15)(112,30)

\put(170,40){\ellipse{100}{35}}%\put(160,60){\oval(100,35)}
\put(118,37){$\bullet$}%\put(120,60){\circle*{5}}
\put(140,23){$\bullet$}%\put(140,46){\circle*{5}}
\put(196,23){$\bullet$}%\put(196,45){\circle*{5}}

\put(140,26){\vector(4,-1){35}}
\put(165,8){$v$}
\put(140,30){$b$}
\put(100,28){\ovalbox{\large $\gamma$}$\leadsto$}

 \end{picture}

         \caption{Curve $\gamma$ and its horizontal lift $\Gamma_{r}$.}
         \label{fig45cex422}
         \index{pictures!Counterexample\ref{lcex422}}
      \end{figure}

It can be shown (see~\cite{CAS}) that the unique horizontal 
lift $\Gamma_{r}$ (see Fig.~{\ref{fig45cex422}}) of $\gamma$
through $\xi_{r}$, where $r>0$, can be parametrized by,
\[
\Gamma_{r}(t)= \left(\gamma_{1}(t)\cdot\sqrt{r^{2}+1},
\gamma_{2}(t)\cdot\sqrt{r^{2}+1},r\right)
\hspace{0.1in} \forall t\in[0,1] %\hspace{0.5in}\ce{2}
\]

%(see Fig.~{\ref{fig45cex422}})

Notice that $M\setminus F_{b}\neq \emptyset$.

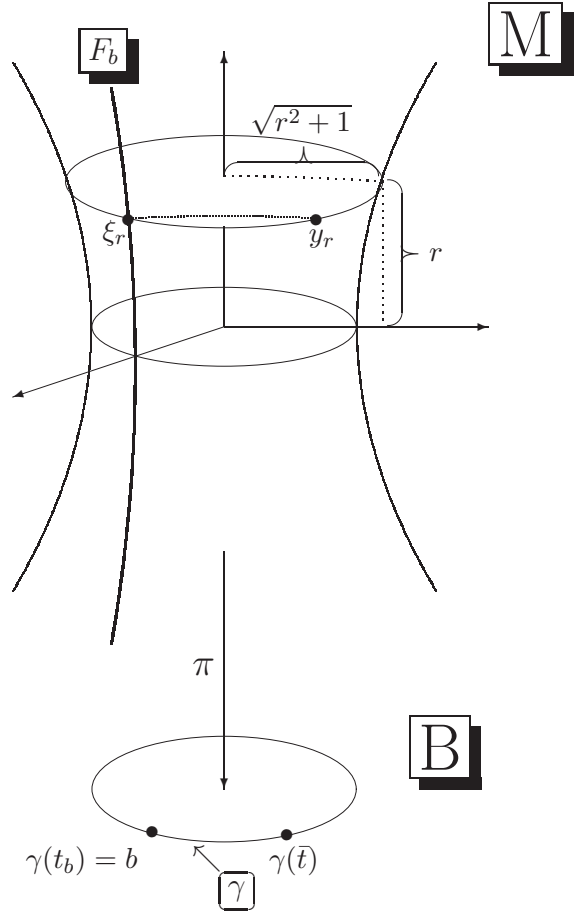
\begin{figure}[here]

         %\begin{picture}(100,300) (-60,10)
           %\includegraphics{piccex422-6.ps}
        \begin{picture}(305,340)(0,0)

%\dottedline{2}(0,350)(305,350)
%\dottedline{2}(0,0)(305,0)

\put(270,320){\shadowbox{\Huge M}}
%parallels=horizontal curves

%\qbezier[90](110,300)(170,320)(230,300)
%\qbezier(110,300)(170,280)(230,300)
\put(170,280){\ellipse{120}{35}}

%\qbezier[90](120,245)(170,260)(220,245)
%\qbezier(120,245)(170,230)(220,245)
\put(170,225){\ellipse{100}{30}}

%meridians=vertical curves
\qbezier(90,125)(149,225)(90,325)

\put(115,317){\shadowbox{\large $F_{b}$}}
\thicklines
\qbezier(127,105)(145,210)(127,315)
\thinlines

\qbezier(250,125)(190,225)(250,325)
%axis
\put(170,225){\line(0,1){38}}
\put(170,282){\vector(0,1){47}}
\put(170,225){\vector(1,0){100}}
\put(170,225){\vector(-3,-1){80}}
%text
\put(131,263){$\bullet$}%\put(134,286){\circle*{5}}
\put(202,263){$\bullet$}%\put(205,286){\circle*{5}}

\put(124,258){$\xi_{r}$}
\put(202,257){$y_{r}$}
\qbezier[60](134,266)(171,268)(205,266)
%\dashline[+90]{3}(143,285)(200,285)

\put(180,299){$\sqrt{r^{2}+1}$}
\put(196,287){$\curlywedge$}
\qbezier[20](170,282)(200,282)(230,280)
\put(199,282){\oval(58,10)[t]}

\put(236,250){$\succ r$}%\put(237,270){$\leadsto r$}
\qbezier[20](230,280)(230,252)(230,225)
%\put(230,271){\vector(0,1){30}}
\put(232,253){\oval(10,55)[r]}

\put(240,50){\shadowbox{\Huge B}}

\put(170,140){\vector(0,-1){90}}
\put(158,93){\Large $\pi$}

%\qbezier(38,30)(75,45)(112,30)
%\qbezier(38,30)(75,15)(112,30)

%\put(170,60){\oval(100,35)}
\put(170,50){\ellipse{100}{40}}
\put(140,31){$\bullet$}%\put(140,44){\circle*{5}}
\put(191,30){$\bullet$}%\put(191,43){\circle*{5}}

\put(95,20){$\gamma(t_{b})=b$}
\put(168,10){\ovalbox{\large $\gamma$}}
\put(158,21){$\nwarrow$}
\put(187,20){$\gamma(\bar{t})$}

 \end{picture}

         \caption{Generic elements $\xi_{r}$ in the fiber $F_{b}$,
                  and $y_{r}$ in $M\setminus F_{b}$.}
         \label{fig46cex422}
         \index{pictures!Counterexample\ref{lcex422}}
      \end{figure}

Now, any element $y_{r}$ of $M\setminus F_{b}$ is of the form,
\[
y_{r}=\left(\sqrt{r^{2}+1}\cdot\cos\bar{t},
            \sqrt{r^{2}+1}\cdot\sin\bar{t},r\right)
\]
for some $r\in\real$ and $\bar{t}\in [0,2\pi%](
), \bar{t}\neq t_{b}$, where $\gamma(t_{b})=(\cos t_{b},\sin t_{b},0)=
b\neq \pi(y_{r})=(\cos\bar{t},\sin\bar{t},0)=\gamma(\bar{t})$
(see Fig.~{\ref{fig46cex422}}).

%(see Fig.~{\ref{fig46cex422}})

\vspace{0.1in}

We may choose $\bar{t}\in [0,2\pi%](
)$ as follows,
\[
\bar{t}:=
\left\{
\begin{array}{cccc}
t_{b}+\pi, & \mbox{ if }& 0\leq t_{b}<\pi & (:. \pi\leq \bar{t}<2\pi)\\
t_{b}-\pi, & \mbox{ if }& \pi\leq t_{b}<2\pi & (:. 0\leq \bar{t}<\pi)
\end{array}
\right.%}
\]

In particular, $\bar{t}\neq t_{b}$ and $|\bar{t}-t_{b}|=\pi$.

Moreover,
\begin{eqnarray}
\protect\label{leqn441}
\lefteqn{d_{M}(\xi_{r}, y_{r})
=\sqrt{\left(r^{2}+1\right)(\cos t_{b}-\cos\bar{t})^{2}+
\left(r^{2}+1\right)(\sin t_{b}-\sin\bar{t})^{2}}=}\nonumber\\
&=&\sqrt{r^{2}+1}
\sqrt{1-2\cos t_{b}\cos\bar{t}-2\sin t_{b}\sin\bar{t}+1}=
\nonumber\\
&=&\sqrt{r^{2}+1}\sqrt{2}\sqrt{1-\cos(t_{b}-\bar{t})}=
\sqrt{r^{2}+1}\sqrt{2}\sqrt{2}=2\sqrt{r^{2}+1}
			  %\hspace{0.2in}\mbox{\ce{7}}
\end{eqnarray}

Let $\epsilon >0$ be arbitrary.

\vspace{0.1in}

\begin{figure}[here]

       \begin{picture}(360,330)(-10,15)

%\dottedline{2}(0,355)(360,355)
%\dottedline{2}(0,10)(360,10)

\put(315,319){\shadowbox{\huge M}}
\dashline[+90]{3}(0,350)(350,350)
\dashline[+90]{3}(0,98)(0,350)
\dashline[+90]{3}(350,98)(350,350)
\dashline[+90]{3}(0,98)(350,98)

%parallels=horizontal curves

\dottedline{2}(88,270)(257,290)
\dottedline{2}(88,270)(254,270)

%meridians=vertical curves
\thicklines
\qbezier(90,105)(100,340)(50,345)
\put(77,316){\shadowbox{\Large $F_{b}$}}

\thinlines
\qbezier(250,105)(240,340)(290,345)
%axis
\put(170,105){\vector(0,1){235}}
       \put(168,343){$e_{3}$}
\put(15,135){\vector(1,0){305}}
       \put(325,134){$e_{2}$}
%text
\put(84,267){$\bullet$}%\put(88,270){\circle*{5}}
   \put(75,266){$\xi_{r}$}
   %\put(87,269){\begin{rotate}{-40}
   %              $\Box$
   %              \end{rotate}}
   %\put(90,270){$\cdot$}
   %\put(87,267){$\boxdot$}
   \thicklines
   \put(86,278){\line(1,-4){2}}
   \put(94,280){\line(-4,-1){8}}
   \put(96,272){\line(-1,4){2}}
   \put(88,270){\line(4,1){8}}
   \put(89,272){$\cdot$}

\put(250,267){$\bullet$}%\put(254,270){\circle*{5}}
   \put(260,268){$y_{r}$}
\put(254,287){$\bullet$}%\put(257,290){\circle*{5}}
   \put(262,288){$y_{\epsilon}(\epsilon,b)$}
\put(89,132){$\bullet$}%\put(91,135){\circle*{5}}
   \put(3,153.5){\ovalbox{$(\cos t_{b}, \sin t_{b}, 0)$}}
   \put(78,139){$\searrow$}
   \put(91,135){$\boxdot$}
%  \put(91,135){$\Box$}\put(93,135){$\cdot$}
\put(246,132){$\bullet$}%\put(249,135){\circle*{5}}
\put(254,153.5){\ovalbox{$(\cos \bar{t}, \sin \bar{t}, 0)$}}
   \put(251,139){$\swarrow$}
   \put(241,135){$\boxdot$}
%  \put(241,135){$\Box$}\put(244,135){$\cdot$}
\put(179,200){$\succ |r|$}%\put(182,200){$\leadsto |r|$}
\put(170,202){\oval(20,135)[r]}

\put(131,134){\oval(76,20)[b]}
\put(124,118){$\curlyvee$}
\put(127,107){1}

\put(208,134){\oval(76,20)[b]}
\put(205,118){$\curlyvee$}
\put(207,107){1}

\put(295,50){\shadowbox{\huge B}}

\put(170,89){\vector(0,-1){40}}
\put(160,71){\large $\pi$}

\put(169,46){\ellipse{160}{40}}
\put(87,43){$\bullet$}%\put(89,46){\circle*{5}}
\put(245,44){$\bullet$}%\put(249,46){\circle*{5}}

\put(37,46){$\gamma(t_{b})=b$}
\put(78,25){\ovalbox{\large $\gamma$}}
\put(91,31){$\leadsto$}
%\put(156,32){$\searrow$}
\put(255,46){$\gamma(\bar{t})$}

         \end{picture}

         \caption{Realization of the distance between $y_{\epsilon}$
                  and the fiber $F_{b}$.}
         \label{fig47cex422}
         \index{pictures!Counterexample\ref{lcex422}}
      \end{figure}
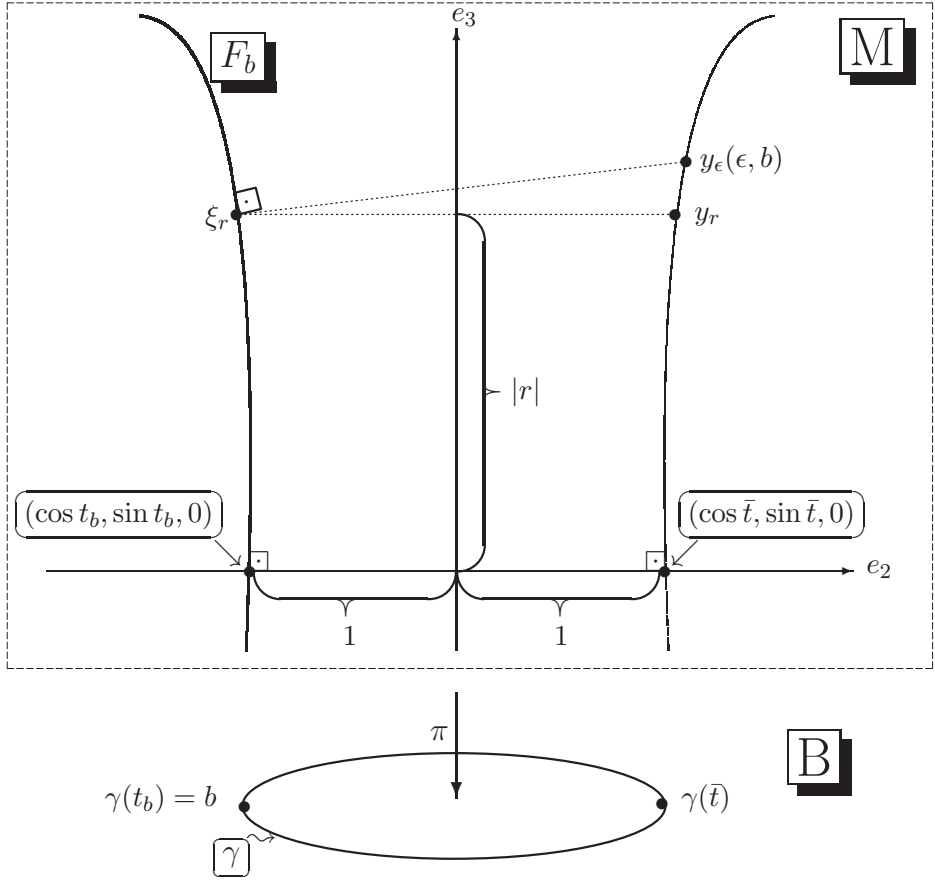

If $\epsilon\leq 2$, by (\ref{leqn441}) we have
(see Fig.~{\ref{fig47cex422}}),
\[
d_{M}(\iota(F_{b}),y_{0}) = d_{M}(\xi_{0}, y_{0})
\stackrel{(\ref{leqn441})}{=} 2 \geq \epsilon
\]
which shows that {\bf [RI.2]} fails for $y_{\epsilon}:=y_{0}=
(\cos\bar{t}, \sin\bar{t}, 0)\in M\setminus F_{b}$.

If $\epsilon>2$, consider any $r\in \real$ satisfying
\[
r > \displaystyle{\frac{\sqrt{\epsilon^{2}-4}}{2}}=
 \displaystyle{\frac{\sqrt{\epsilon-2}\sqrt{\epsilon+2}}{2}}>0
\]
this choice of $r$ being possible, because property~(\ref{leqn440}) 
holds for this Counterexample.

In that case, 
we have,
\begin{eqnarray}
\protect\label{leqn442} 
2r>\sqrt{\epsilon^{2}-4}>0 & \Longrightarrow & 
4r^{2}>\epsilon^{2}-4>0
\Longrightarrow 4(r^{2}+1)>\epsilon^{2}\Longrightarrow\nonumber\\
&  \Longrightarrow  & 2\sqrt{r^{2}+1}>\epsilon
%\hspace{0.5in}\mbox{\ce{8}}
\end{eqnarray}

In what follows we will define (see Fig.~{\ref{fig47cex422}}),
\[
y_{\epsilon}=\left(\sqrt{r_{\epsilon}^{2}+1}\cdot\cos\bar{t},
\sqrt{r_{\epsilon}^{2}+1}\cdot\sin\bar{t},r_{\epsilon}\right)\in
M\setminus F_{b} %\hspace{0.5in}{\bf [CE.??]}
\]
satisfying the 2 conditions,
\begin{itemize}
\item $r_{\epsilon}>r$; and
\item the unique straight line passing through $y_{\epsilon}$ and
$\xi_{r}$ is perpendicular to $F_{b}$ at $\xi_{r}$, thus giving us the
realization of the distance $d_{M}(y_{\epsilon},\iota(F_{b}))=
d_{M}(y_{\epsilon,\xi_{r}})$.
\end{itemize}

%(see Fig.~{\ref{fig47cex422}})

We may assume for the sake of a much simplified calculation,
that, $t_{b}=\displaystyle{\frac{3\pi}{2}}$ and $\bar{t}=
\displaystyle{\frac{\pi}{2}}$, since $M$ is
symmetric with respect to both axis $e_{3}$ and $e_{2}$

Thus we have, $\gamma\left(\frac{3\pi}{2}\right)=
(0,-1,0)=b, \gamma\left(\frac{\pi}{2}\right)=
(0,1,0)$, a generic element $\xi_{r}=(0, -\sqrt{r^{2}+1}, r)$ in the
fiber $F_{b}$, and a generic element $y_{r}=(0, \sqrt{r^{2}+1}, r)$
in $M$, but not in the fiber $F_{b}$.

In this case, $F_{b}$ is given by $x_{2}= -\sqrt{x_{3}^{2}+1}$, and the
perpendicular line to $F_{b}$ at $\xi_{r}$ has equation,
\begin{eqnarray}
\protect\label{leqn443}
x_{2}+\sqrt{r^{2}+1}&=&\displaystyle{\frac{-1}{\mu}}(x_{3}-r)=
\nonumber\\
&=&\displaystyle{\frac{\sqrt{r^{2}+1}}{r}}(x_{3}-r),
\hspace{0.2in}\forall x_{3}\in \real %\hspace{0.7in}\mbox{\ce{9}}
\end{eqnarray}
where $\mu=\displaystyle{\frac{d}{du}} -\sqrt{u^{2}+1}\mid_{u=r}=
\displaystyle{\frac{-2r}{2\sqrt{r^{2}+1}}}$.

Since $y_{\epsilon}=\left(0, \sqrt{r_{\epsilon}^{2}+1}, r_{\epsilon}
\right)$ is on the line (\ref{leqn443}), we obtain the following
equation,
\begin{eqnarray}
\protect\label{leqn444}
\sqrt{r_{\epsilon}^{2}+1}+\sqrt{r^{2}+1} & = & 
\displaystyle{\frac{\sqrt{r^{2}+1}}{r}}(r_{\epsilon}-r)>0
 %\hspace{0.5in}\mbox{\ce{10}}
\end{eqnarray}
which defines $r_{\epsilon}$.

Indeed, equation (\ref{leqn444}) has only one solution,
\[
\sqrt{r_{\epsilon}^{2}+1}=\displaystyle{\frac{\sqrt{r^{2}+1}}{r}}
r_{\epsilon}-2\sqrt{r^{2}+1}=\sqrt{r^{2}+1}\left(
\displaystyle{\frac{r_{\epsilon}}{r}}-2\right)\Rightarrow
\]
\[
\Rightarrow  r_{\epsilon}^{2}+1=(r^{2}+1)\left(\displaystyle{
\frac{r_{\epsilon}^{2}}{r^{2}}}+4-4\displaystyle{
\frac{r_{\epsilon}}{r}}\right) \Rightarrow
\]
\[
\Rightarrow r_{\epsilon}^{2}+1=r_{\epsilon}^{2}+\displaystyle{
\frac{r_{\epsilon}^{2}}{r^{2}}}+4(r^{2}+1)-4\left(r+\frac{1}{r}\right)
r_{\epsilon}\Rightarrow
\]
\[
\Rightarrow \displaystyle{\frac{1}{r^{2}}}{\bf r_{\epsilon}^{2}}
-4\left(r+\frac{1}{r}\right){\bf r_{\epsilon}}+4r^{2}+3=0\Rightarrow
\]
\[
\Rightarrow {\bf r_{\epsilon}}=\displaystyle{\frac{4\left(r+
\displaystyle{\frac{1}{r}}\right)
\pm\sqrt{\Delta}}{\displaystyle{\frac{2}{r^{2}}}}}
\]
where
\begin{eqnarray*}
\Delta&=&16\left(r+\displaystyle{\frac{1}{r}}\right)^{2}
-4\displaystyle{\frac{1}{r^{2}}}\left(4r^{2}+3\right)=
16\left(r^{2}+2+\displaystyle{\frac{1}{r^{2}}}\right)-16-\displaystyle{
\frac{12}{r^{2}}}=\\
&=& 16r^{2}+16+\displaystyle{\frac{4}{r^{2}}}=
4\left(4r^{2}+4+\displaystyle{\frac{1}{r^{2}}}\right)=
4\left(2r+\displaystyle{\frac{1}{r}}\right)^{2}\Rightarrow\\
 \Rightarrow
\sqrt{\Delta}&=&2\left(2r+\displaystyle{\frac{1}{r}}\right)=
\left(4r+\displaystyle{\frac{2}{r}}\right)
\end{eqnarray*}
which implies that,
\[
{\bf r_{\epsilon}}=\displaystyle{
\frac{4\left(r+\displaystyle{\frac{1}{r}}\right)
\pm
\left(4r+\displaystyle{\frac{2}{r}}\right)}
{\displaystyle{\frac{2}{r^{2}}}}}=2r^{3}+2r\pm\left(2r^{3}+r\right)=
\left\{
\begin{array}{l}
4^{3}+3r  \mbox{ (solution) }\\
r  \mbox{ (\underline{not} a solution) }
\end{array}
\right.%}
\]
Then, the only solution of equation (\ref{leqn444}) is,
\begin{eqnarray}
\protect\label{leqn445}
{\bf r_{\epsilon}}=4^{3}+3r=r\left(4^{2}+3\right)>r>0
 %\hspace{0.5in}\mbox{\ce{11}}
\end{eqnarray}

Finally, by employing the expressions,
\[
y_{\epsilon}=\left(0,
\sqrt{r_{\epsilon}^{2}+1}, r_{\epsilon}\right)
\hspace{0.1in}\mbox{ and }\hspace{0.1in}
\xi_{r}=(0, -\sqrt{r^{2}+1}, r)
\]and by using (\ref{leqn444}) and (\ref{leqn445}) , we can
estimate the distance from $y_{\epsilon}$ to the fiber $F_{b}$, only
in terms of $r$,
\begin{eqnarray}
\protect\label{leqn446}
\lefteqn{d_{M}(y_{\epsilon},\xi_{r})=\sqrt{0^{2}+
\left(\sqrt{r_{\epsilon}^{2}+1}+\sqrt{r^{2}+1}\right)^{2}
+(r_{\epsilon}-r)^{2}} \stackrel{(\ref{leqn444})}{=}}\nonumber\\
&=&\sqrt{\left(\displaystyle{\frac{r^{2}+1}{r^{2}}}\right)
(r_{\epsilon}-r)^{2}+(r_{\epsilon}-r)^{2}}=
\sqrt{\left(\displaystyle{\frac{r^{2}+1}{r^{2}}}+1\right)
(r_{\epsilon}-r)^{2}}=\nonumber\\
&=&\sqrt{\displaystyle{\frac{\left(2r^{2}+1\right)}{r^{2}}}
(r_{\epsilon}-r)^{2}}\stackrel{(\ref{leqn445})}{=}
\sqrt{\displaystyle{\frac{\left(2r^{2}+1\right)}{r^{2}}}
(4r^{3}+3r-r)^{2}}=\nonumber\\
&=&\sqrt{\left(2r^{2}+1\right)\left(4r^{2}+2\right)}=
\sqrt{2r^{2}+1}\left(4r^{2}+2\right)=2\sqrt{2r^{2}+1}
\left(2r^{2}+1\right)=\nonumber\\
&=& 2\left(2r^{2}+1\right)^{\frac{3}{2}} 
%\hspace{2.0in}\mbox{\ce{12}}
\end{eqnarray}

Next, we claim that,
\begin{eqnarray}
\protect\label{leqn447}
\left(2r^{2}+1\right)^{\frac{3}{2}}>\sqrt{r^{2}+1}
%\hspace{0.5in}\mbox{\ce{13}}
\end{eqnarray}

The function  $f\in\cal{C}^{\infty}(\real)$ defined by,
\[
f(r):=\left(2r^{2}+1\right)^{3}-\left(r^{2}+1\right)
\]
is clearly strictly increasing on $[0,%](
\infty)$ 
, and thus,
\[
f(r)>f(0)=0, \hspace{0.1in} \forall r>0 \Rightarrow
\left(2r^{2}+1\right)^{\frac{3}{2}}>
\left(r^{2}+1\right)^{\frac{1}{2}}, \hspace{0.1in} \forall r>0
\]
which is claim (\ref{leqn447}).

Now, if we combine (\ref{leqn446}), (\ref{leqn447}) and 
(\ref{leqn442}),
we get,
\[
d_{M}(y_{\epsilon},\iota(F_{b}))=d_{M}(y_{\epsilon},\xi_{r})>2
\sqrt{r^{2}+1}>\epsilon
\]
which shows that {\bf [RI.2]} fails for $y_{\epsilon}:=
\left(0,
\sqrt{r_{\epsilon}^{2}+1}, r_{\epsilon}\right)
\in M\setminus F_{b}$.

We have thus shown that for any $\epsilon >0$ 
there exists $y_{\epsilon}\in M\setminus F_{b}$ 
for which {\bf [RI.2]} fails. Consequently, the inclusion 
map $\iota:F_{b}\rightarrow M$ is \underline{not} a rough isometry.

This describes the Counterexample.

\vspace{0.1in}

Next, including a lower bound in assumption (\ref{leqn416}), and adding 
an universal diameter upper bound condition on the fibers, we will 
show that $\pi: M\rightarrow B$ is a rough isometry.

\begin{thm423}
\protect\label{lthm423}
 Let $\pi: M\rightarrow B$ be an onto smooth map with
maximal rank, where $B$ is complete. Assume the following,
\begin{description}
\item[(UDF)]   $\exists m>0$, a universal constant, such that 
$\mbox{diam }\left\{\pi^{-1}(b)\right\} \leq m < \infty$,
for all $b\in B$; and 
\item[(HLC)] $\exists \alpha\geq 1$ and $\beta >0$ such that,  for all 
$b\in B$ the inequality holds:
\[
\frac{1}{\alpha} ||w||_{B} - \beta \leq ||v||_{M} \leq \alpha ||w||_{B}
+ \beta
\]
for all $x\in F_{b}$ and $w\in T_{b}B$,
where $v\in (HT)_{x}\subset T_{x}M$
is the horizontal lift of $w$ through $x$ and we assume that
$v$ satisfies $||v||_{M}\leq 1$.
\end{description}
Then, $\pi: M\rightarrow B$ is a rough isometry.

In particular, if the fiber $\pi^{-1}(b_{0})$ is compact for 
some $b_{0}$, then $M$ is roughly isometric to the product 
$\pi^{-1}(b_{0})\times B$.
\end{thm423}
\pf\hspace{0.1in}
Firstly,  note that in \textbf{(HLC)} the horizontal lift $v\in (HT)_{x}$ 
of $w$ is assumed to satisfy $||v||_{M}\leq 1$.
Otherwise, if $||v||_{M}>1$ we define
$\tilde{v}:=\displaystyle{\frac{v}{||v||_{M}}}$, with the properties
\begin{itemize}
\item $\tilde{v}:=\displaystyle{\frac{v}{||v||_{M}}\in (HT)_{x}}$
\item $||\tilde{v}||_{M}=1$
\item $(\pi_{\ast})_{x}(\tilde{v})=\displaystyle{\frac{w}{||v||_{M}}}$
\end{itemize}
and if we use $\displaystyle{\frac{w}{||v||_{M}}}$ and $\tilde{v}$
in \textbf{(HLC)}, we thus obtain the equivalent inequality,
\begin{eqnarray*}
& &
\frac{1}{\alpha}\left|\left|\frac{w}{||v||_{M}}\right|\right|_{B}-\beta
\leq ||\tilde{v}||_{M}\leq
\alpha\left|\left|\frac{w}{||v||_{M}}\right|\right|_{B}+\beta
\Rightarrow\\
& \Rightarrow &
\frac{1}{\alpha} \frac{||w||_{B}}{||v||_{M}} - \beta
\leq \frac{||v||_{M}}{||v||_{M}}\leq
\alpha\frac{||w||_{B}}{||v||_{M}} + \beta
\Rightarrow \\
& \Rightarrow &
\frac{1}{\alpha} ||w||_{B} - \beta||v||_{M} \leq ||v||_{M}\leq
\alpha ||w||_{B} + \beta||v||_{M}
\Rightarrow \\
& \Rightarrow &
\frac{1}{\alpha} ||w||_{B} \leq \left(\beta +1\right)||v||_{M}\wedge
\left(1-\beta\right)||v||_{M}\leq \alpha ||w||_{B}
\Rightarrow\\
& \Rightarrow &
\left\{
\begin{array}{cl}
\displaystyle{\frac{1}{\alpha\left(\beta +1\right)}} ||w||_{B}
\leq ||v||_{M}\leq
\displaystyle{\frac{\alpha}{\left(1-\beta\right)}} ||w||_{B},
& \mbox{ if }\beta\neq 1\\
\displaystyle{\frac{1}{\alpha\left(\beta +1\right)}} ||w||_{B}
\leq ||v||_{M},
& \mbox{ if }\beta = 1
\end{array}
\right.%}
\end{eqnarray*}
for $w\in T_{b}B$,  where $v$ is the unique horizontal lift of
$w$ through $x$ with $||v||_{M}>1$.

We must verify the validity of {\bf (RI.1)} and
{\bf (RI.2)}.

Clearly, axiom {\bf (RI.2)} holds since $\pi$ is onto.

To verify {\bf (RI.1)}, let $x,y\in M$.

We may assume that $B$ is connected. 
Otherwise, we repeat the argument which will be utilized in this 
proof, on each connected component and the result will follow.

Because $B$ is complete, there exists a minimal geodesic
$\gamma$ joining $\pi(x)$ to $\pi(y)$, with 
$\ell(\gamma)=d_{B}(\pi(x),\pi(y))$, which we parametrize by
$\gamma:[0,1]\rightarrow~B,$ where,
$\gamma(0) := \pi(x), \gamma(1) := \pi(y)$.

Recall that $\gamma$ has a unique horizontal lift 
$\Gamma_{x}:[0,1]\rightarrow M$, through $x$, so $\Gamma_{x}$ 
intersects the fiber $F_{\pi(y)}$ containing $y$.

We may assume, without loss of generality, that $\Gamma_{x}$ is 
parametrized proportionally to arc lenght and
$||{\Gamma_{x}}^{\prime}(t)||_{M}\leq 1$ for all $t\in [0,1]$.

Thus we can write,
\begin{eqnarray}
\protect\label{leqn420}
\ell(\Gamma_{x}) & = & \int_{0}^{1}||\Gamma'_{x}||_{M} dt
\stackrel{(HLC)} 
 \leq  \alpha \int_{0}^{1}
||(\pi_{\ast})_{\Gamma_{x}(t)}\Gamma'_{x}(t)||_{B}  dt
+ \beta
\stackrel{\pi\circ\Gamma_{x}}{=}  \nonumber \\
& = & \alpha \int_{0}^{1} ||\gamma'(t)||_{B}  dt + \beta
=
\alpha\cdot \ell(\gamma) + \beta
\stackrel{\mbox{\footnotesize  def }\gamma}{=} \nonumber\\
& = & \alpha\cdot d_{B}(\pi(x),\pi(y)) + \beta
\end{eqnarray}

By the triangle inequality, by hypothesis and the above, we have,
\begin{eqnarray*}
d_{M}(x,y) & \stackrel{\triangle}{\leq} & 
d_{M}(x,\Gamma_{x}(1)) + d_{M}(\Gamma_{x}(1),y) 
\stackrel{\mbox{\footnotesize dist.}}{\leq} \\
             & \leq & \ell(\Gamma_{x}) + d_{M}(\Gamma_{x}(1),y)
             \stackrel{\bf (UDF)}{\leq} \ell(\Gamma_{x}) + m
             \stackrel{(\ref{leqn420})}{\leq} \\
             & \leq & \alpha\cdot d_{B}(\pi(x),\pi(y)) + \beta
             + m
\end{eqnarray*}
which can be rewritten as,
\begin{equation}
\protect\label{lequ421}
d_{B}(\pi(x),\pi(y)) \geq \frac{1}{\alpha}
d_{M}(x,y) - 
\frac{(\beta + m)}{\alpha}
\end{equation}

Now, we claim that for $\gamma$, the minimal geodesic joining
$\pi(x)$ to $\pi(y)$, its length $\ell(\gamma)$ satisfies,
\begin{equation}
\protect\label{lequ422}
d_{B}(\pi(x),\pi(y)) = \ell(\gamma) \leq
\alpha\cdot \ell(\varsigma) + \alpha\cdot\beta
\end{equation}
for any smooth curve $\varsigma: [0,1]\rightarrow M$, joining $x$
to $y$.

First, observe that for any orthogonal vectors $U$ and $W$,
\[
||U\oplus W||^{2}=||U||^{2} + ||W||^{2}\geq
\max\{||U||^{2},||W||^{2}\}
\]

Now, since each tangent vector is the direct sum of a horizontal
and a vertical vector, we can write,
\begin{equation}
\protect\label{lequ423}
 \ell(\varsigma) = \int_{0}^{1}
||\varsigma'(t)||_{M}  dt = \int_{0}^{1} ||\varsigma'_{H}(t)\oplus
\varsigma'_{V}(t)||_{M} dt \geq \int_{0}^{1}
||\varsigma'_{H}(t)||_{M}  dt
\end{equation}
where we are assuming here that $\varsigma_{H}$ is parametrized 
proportional to arclength, and
$||\varsigma'_{H}(t)||_{M} \leq 1, \forall t\in [0,1]$.

Since, $(VT)_{x} = \ker (\pi_{\ast})_{x}, \forall x\in M$, we have,
\begin{equation}
\protect\label{lequ424}
\ell(\pi\circ\varsigma) = \int_{0}^{1}
||(\pi_{\ast})_{\varsigma(t)}\varsigma'(t)||_{B} dt = \int_{0}^{1}
||(\pi_{\ast})_{\varsigma(t)}\varsigma'_{H}(t)||_{B} dt
\end{equation}

From the left-hand side of {\bf (HLC)},
\begin{eqnarray}
\protect\label{leqn425}
 & & \frac{1}{\alpha}\cdot ||(\pi_{\ast})_{\varsigma(t)}
 \varsigma'_{H}(t)||_{B} - \beta
   \leq ||\varsigma'_{H}(t)||_{M} \Rightarrow \nonumber\\
 & \Rightarrow &
 ||(\pi_{\ast})_{\varsigma(t)}\varsigma'_{H}(t)||_{B}\leq
 \alpha\cdot
 ||\varsigma'_{H}(t)||_{M} + \alpha\cdot\beta
\end{eqnarray}
for all $t\in [0,1]$.

If we combine (\ref{lequ423}), (\ref{lequ424}) and
(\ref{leqn425}), we get,
\begin{eqnarray}
\protect\label{leqn427}
\ell(\pi\circ\varsigma) & \stackrel{(\ref{lequ424})}{=} &
\int_{0}^{1}||(\pi_{\ast})_{\varsigma(t)}\varsigma'_{H}(t)||_{B}
dt \stackrel{(\ref{leqn425})}{\leq}
\alpha\int_{0}^{1}||\varsigma'_{H}(t)||_{M} + \alpha\beta \leq 
\nonumber \\
& \stackrel{(\ref{lequ423})}{\leq} & \alpha \ell(\varsigma) +
\alpha\beta
\end{eqnarray}

Inequality (\ref{leqn427}) and the fact that $\gamma$ is a minimal
geodesic joining $\pi(x)$ to $\pi(y)$ imply that we can finally
write,
\[
d_{B}(\pi(x),\pi(y)) = \ell(\gamma) \leq
\ell(\pi\circ\varsigma)\leq\alpha \ell(\varsigma) + \alpha\beta
\]
for any smooth curve $\varsigma: [0,1]\rightarrow M$, joining $x$
to $y$, which is claim (\ref{lequ422}).

We recall that by definition of infimum, $d_{M}(x,y)$ is the
greatest lower bound for
 $\{\ell(\varsigma), \mbox{ where }
\varsigma:[0,1]\rightarrow M \mbox{ is any smooth curve joining }
x \mbox{ to }y\}$, and since $\varsigma$ is arbitrary in
(\ref{lequ422}), we obtain,
\begin{equation}
\protect\label{lequ426}
d_{B}(\pi(x),\pi(y)) \leq \alpha\cdot
d_{M}(x,y) + \alpha\cdot\beta
\end{equation}

Let $A:=\alpha\geq 1$ and $C:=\max\left\{
\displaystyle{\frac{\beta + m}{\alpha}},\alpha\cdot\beta\right\}>0$.

If we now, rewrite (\ref{lequ421}) and (\ref{lequ426}) in terms of
$A$ and $C$, as follows,
\begin{eqnarray*}
\frac{1}{A} d_{M}(x,y) - C & \leq & \frac{1}{\alpha} d_{M}(x,y) -
\frac{(\beta + m)}{\alpha}
\stackrel{(\ref{lequ421})}{\leq} \\
& \leq & d_{B}(\pi(x),\pi(y)) \stackrel{(\ref{lequ426})}{\leq}
\alpha d_{M}(x,y) + \alpha\beta \leq A d_{M}(x,y) + C
\end{eqnarray*}
we obtain {\bf (RI.1)} for $\pi$.

\pfe

\vspace{0.1in}

In the following two Counterexamples,  we show that the universal 
diameter property \textbf{(UDF)} of the fibers, and the control over 
the length of horizontal lifts \textbf{(HLC)} of tangent vectors 
are both necessary conditions in {\bf Theorem~\ref{lthm423}}. 

\begin{cex424}
\protect\label{lcex424}
We will exhibit $M,B,\pi$, where $B$ is connected and each fiber $F_{b}$ 
is compact for all $b\in B$, satisfying all but condition {\bf (HLC)} 
in \textsf{Theorem~\ref{lthm423}},
i.e.,

For any given constants $\alpha\geq 1$ and $\beta >0$, there exist
$\bar{b}\in~B$, $\bar{x}\in~F_{\bar{b}}$, $\bar{w}\in~T_{\bar{b}}B$ such
that either one of the following holds:
\begin{eqnarray}
\protect\label{leqn448}
||\bar{v}||_{M} > \alpha ||\bar{w}||_{B} + \beta  & \mbox{ or } & 
||\bar{v}||_{M} < \frac{1}{\alpha} ||\bar{w}||_{B} - \beta
%\hspace{0.5in}\mbox{\ce{14}}
\end{eqnarray} 
where $\bar{v}$ is the unique horizontal lift of $\bar{w}$ through
$\bar{x}$.

In this case, the map $\pi: M\rightarrow B$ is not a rough 
isometry.
\end{cex424}

Let $M=\{(x,y,z)\in \real^{3}: x^{2}+z^{2}=1,
y\in \real\}$ and $B=\real$, a \underline{complete} and
\underline{connected} Riemannian manifold.

      \begin{figure}[here]

         %\begin{picture}(100,360) (0,0)
            %\includegraphics{piccex424-9.ps}
         \begin{picture}(355,200)(0,0)

%\dottedline{2}(0,210)(355,210)
%\dottedline{2}(0,0)(355,0)

\put(5,177){\shadowbox{\Huge $M$}}

%parallels=vertical curves
\thicklines
\put(30,135){\ellipse{20}{80}}
\put(70,135){\ellipse{20}{80}}
\put(230,135){\ellipse{20}{80}}
\put(330,135){\ellipse{20}{80}}
%meridians=horizontal curves
\put(30,175){\line(1,0){300}}
\put(30,95){\line(1,0){300}}
\thinlines
\put(236,150){$\bullet$}
\put(242,150){$(x,y,z)$}

\put(80,153){\vector(1,0){34}}
\put(114,149){$X_{\mu}$}
\put(76,150){$\bullet$}
\put(90.5,128.5){\ovalbox{$X(\mu,\eta)$}}
\put(81,142){$\nwarrow$}
\put(80,153){\vector(-1,4){8}}
\put(68,186){$X_{\eta}$}

\put(15,38){\shadowbox{\LARGE $\real$}}

\thicklines
\put(30,35){\line(1,0){300}}
\thinlines
\put(227,33){$\bullet$}
\put(220,40){$y$}

\put(5,6){\shadowbox{\Huge $B$}}

\thicklines
\put(30,3){\line(1,0){300}}
\thinlines
\put(310,0){$\bullet$}
\put(310,10){$f(y)$}

%axis
\put(160,135){\vector(0,1){65}}
       \put(163,193){$e_{3}$}
\put(160,135){\vector(1,0){192}}
       \put(343,140){$e_{2}$}
\put(160,135){\vector(-1,-1){55}}
       \put(116,80){$e_{1}$}

%text
\thicklines
\put(230,90){\vector(0,-1){50}}
\put(235,62){{\large $(x,y,z)\mapsto y$}}
\put(233,32){\vector(3,-1){77}}
\put(247,10){{\Large $f$}}%(295,8)
\put(160,90){\vector(0,-1){85}}
\put(165,43){{\LARGE $\pi$}}
\thinlines
%\qbezier(185,80)(175,45)(185,5)
%\dottedline{2}(88,270)(257,290)

         \end{picture}

      \caption{Manifolds $M,B$ and the map $\pi$ in
      {\bf Counterexample~\ref{lcex424}}.}
         \label{fig49cex424}
         \index{pictures!Counterexample\ref{lcex424}}
      \end{figure}
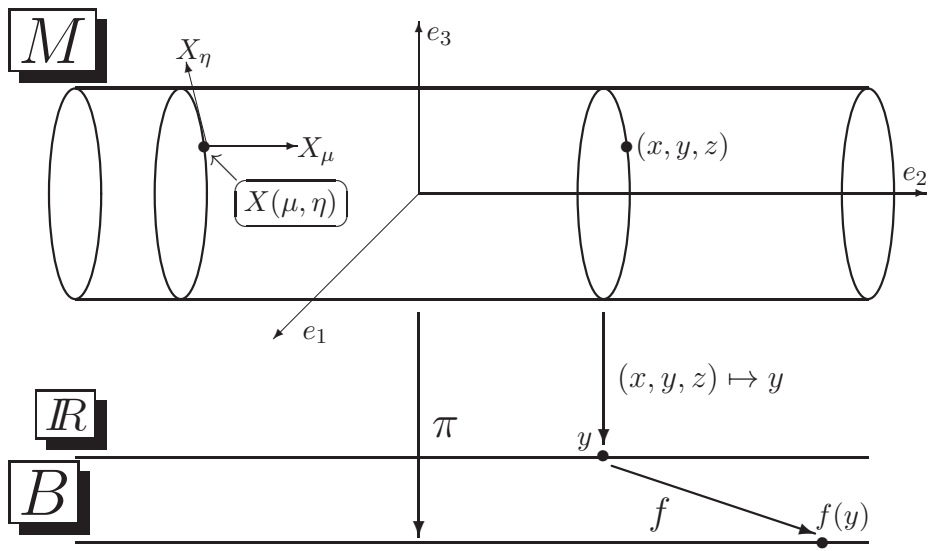

%(see Fig.~{\ref{fig49cex424}})

We first define an auxiliary \ci-diffeomorphism
$f:\real\rightarrow\real$ by,
\begin{eqnarray}
\protect\label{leqn449}
f(y)=\left\{
\begin{array}{lll}
e^{y}-1 & \in [0,\infty%](
) & \mbox{ if }y\geq 0\\
1-e^{-y} & \in (-\infty,0%)[
] & \mbox{ if }y\leq 0
\end{array}\right.%\hspace{0.5in}\mbox{\ce{15}}%\}
\end{eqnarray}

      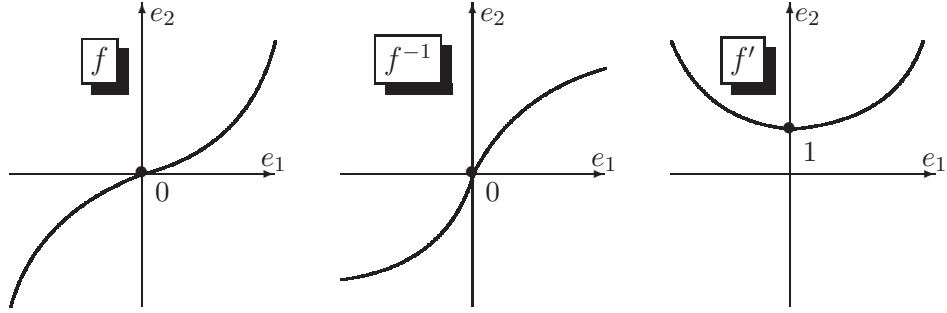
\begin{figure}[here]

         %\begin{picture}(100,360) (0,0)
            %\includegraphics{piccex424-10.ps}
         \begin{picture}(360,100)(0,0)

%\dottedline{2}(0,120)(360,120)
%\dottedline{2}(0,0)(360,0)

%axis
\put(55,0){\vector(0,1){115}}
       \put(58,108){$e_{2}$}
\put(5,50){\vector(1,0){100}}
       \put(100,53){$e_{1}$}
%text f
\thicklines
\qbezier(55,50)(95,60)(105,100)
\qbezier(5,0)(15,35)(55,50)
\thinlines
\put(52,48){$\bullet$}
\put(60,40){0}
\put(32,80){\shadowbox{\large $f$}}%(100,102)

%axis
\put(180,0){\vector(0,1){115}}
       \put(183,108){$e_{2}$}
\put(130,50){\vector(1,0){100}}
       \put(225,53){$e_{1}$}
%text f^{-1}
\thicklines
\qbezier(180,50)(195,80)(230,90)
\qbezier(130,10)(170,15)(180,50)
\thinlines
\put(177,48){$\bullet$}
\put(185,40){0}
\put(143,80){\shadowbox{\large $f^{-1}$}}%(215,100)

%axis
\put(300,0){\vector(0,1){115}}
       \put(302,108){$e_{2}$}
\put(255,50){\vector(1,0){100}}
       \put(350,53){$e_{1}$}
%text f^{\prime}
\thicklines
\qbezier(300,67)(340,70)(350,100)
\qbezier(255,100)(267,70)(300,67)
\thinlines
\put(297,65){$\bullet$}
\put(305,55){1}
\put(274,80){\shadowbox{\large $f^{\prime}$}}

         \end{picture}

      \caption{Graphs of $f,f^{-1},f^{\prime}$.}
         \label{fig410cex424}
         \index{pictures!Counterexample\ref{lcex424}}
      \end{figure}

%(see Fig.~{\ref{fig410cex424}})

Let $\pi:M\rightarrow B$ be given by $\pi(x,y,z):=f(y)$. 

The map $\pi$ is \underline{onto} and \underline{\cin}, since 
both 
the projection $(x,y,z)\mapsto y$ and $f$ have those properties. 
The rank of $\pi$ is \underline{maximal} and no {\bf (HLC)} is 
easily verified (see~\cite{CAS}).

Notice that the fibers have either form,
\[
F_{b}= \pi^{-1}(b)=\left\{ 
\begin{array}{ll}
\left\{(x,\ln(b+1),z)\in \ra{3}: x^{2}+z^{2}=1\right\}, 
& \mbox{ if }b\geq 0\\
\left\{(x,-\ln(1-b),z)\in \ra{3}: x^{2}+z^{2}=1\right\}, 
& \mbox{ if }b < 0
\end{array}
\right.%\}
\] 
Therefore, each fiber $F_{b}$ is \underline{compact} and
$\mbox{diam } F_{b} \leq m,$ for all $b\in B$, where 
$m=3>0$ is the
\underline{universal} upper bound for the fibers' diameters.

\vspace{0.05in}

Finally, we claim that $\pi$ \underline{does not} satisfy 
{\bf (RI.1)}.

\vspace{0.05in}

It suffices to verify that {\bf (RI.1)} fails for $\pi$, for
particular pairs of elements in $M$. We will show that
$\forall A\geq 1, \forall C>0, \exists y_{AC}\in \real$, 
a positive
number such that,
\begin{eqnarray}
\protect\label{leqn4410}
& & d_{B}\left(\pi(x,0,z),\pi(x,y,z)\right)>A\cdot
d_{M}\left((x,0,z),(x,y,z)\right)+C\Leftrightarrow\nonumber\\
& \Leftrightarrow & |f(0)-f(y)|>A\cdot y+C
\Leftrightarrow e^{y}-1>A\cdot y+C %\hspace{0.3in}\ce{18}
\end{eqnarray}
for all $y>y_{AC}$, where $x,z\in\real: x^{2}+z^{2}=1$ are 
arbitrary.

Fixing constants $A\geq 1$ and $C>0$, 
introduce $g\in\cin(\real)$, by
\[
\begin{array}{rcl}
g:y & \mapsto & g(y):=e^{y}-1-Ay-C
\end{array}
\]

One can  show that,
\begin{eqnarray}
\protect\label{leqn4411}
\exists y_{AC}>0: \forall y>y_{AC}\Rightarrow g(y)>0
%\hspace{0.5in}\ce{19}
\end{eqnarray}
using the functional behavior of $g$ (see~\cite{CAS}).

Therefore, (\ref{leqn4410}) holds and the claim follows, and consequently
$\pi$ is \underline{not} a rough isometry.

This describes the Counterexample.

\begin{cex425}
\protect\label{lcex425}
We will exhibit $M,B,\pi$, where $B$ is connected, satisfying all the
conditions in \textsf{Theorem~\ref{lthm423}}, with the exception of
{\bf (UDF)}, 
i.e.,

The fibers' diameters are not uniformly bounded, in other words:
\[
\forall m>0, \exists b_{m} \in B: \mbox{diam } F_{b_{m}} > m
\]

In this case, the map $\pi: M\rightarrow B$ is not a rough isometry.
\end{cex425}
\pf\hspace{0.1in} Let $M=\{(0,y,z)\in\ra{3}\}\cong\{0\}\times\ra{2}$
and $B=\real$, a \underline{complete} and \underline{connected}
Riemannian manifold.

Let $\pi: M\rightarrow B$ be the projection $\pi(x,y,z):=y$. 

The map $\pi$ is 
\underline{onto}, \underline{\cin}, and $\pi$ has
\underline{maximal rank}=1, 
and {\bf (HLC)} is easily verified (see~\cite{CAS}).

\begin{figure}[here]

         %\begin{picture}(100,360) (0,0)
            %\includegraphics{piccex425-12.ps}
         \begin{picture}(350,195)(0,0)

%\dottedline{2}(0,195)(350,195)
%\dottedline{2}(0,0)(350,0)

\put(15,166){\shadowbox{\huge $M$}}
%\dottedline{2}(0,195)(350,195)
%parallels=vertical curves
\put(300,85){\line(0,1){95}}
\thicklines
\put(250,85){\line(0,1){95}}
\thinlines
\put(60,85){\line(0,1){95}}
%meridians=horizontal curves
\put(60,180){\line(1,0){240}}
\thicklines
\put(60,145){\line(1,0){240}}
\thinlines
\put(60,85){\line(1,0){57}}
\put(123,85){\line(1,0){177}}
%\text
\put(157,143){$\bullet$}
\put(121,150){$(0,0,z)$}%(120,162)
%%\put(147,150){$\searrow$}

\put(247,123){$\bullet$}
\put(252,115){$(0,b,0)$}%(263,110)
%%\put(252,115){$\nwarrow$}

\put(247,143){$\bullet$}
\put(211,150){$(0,b,z)$}%(209,162)
%%\put(237,150){$\searrow$}
\put(215,93){\ovalbox{\Large $F_{b}$}$\leadsto$}
\put(250,148){\vector(1,0){30}}
\put(264,152){$w$}
\put(90,130){\vector(1,0){30}}
\put(120,126){$X_{\mu}$}
\put(87,128){$\bullet$}
\put(64,119){$X(\mu,\eta)$}%(70,106)
%%\put(78,118){$\nearrow$}
\put(90,130){\vector(0,1){30}}
\put(85,161){$X_{\eta}$}
%axis
\put(160,125){\vector(0,1){68}}
       \put(163,185){$e_{3}$}
\put(160,125){\vector(1,0){165}}
       \put(320,130){$e_{2}$}
\put(160,125){\vector(-1,-1){55}}
       \put(114,70){$e_{1}$}

\put(15,22){\shadowbox{\huge $B$}}

%manifold
\thicklines
\put(30,17){\line(1,0){300}}
\thinlines
%text
\put(160,82){\vector(0,-1){60}}
\put(165,43){\large{$(x,y,z)\stackrel{\pi}{\mapsto}y$}}
\put(158,15){$\bullet$}
\put(158,5){$0$}
\dashline[+90]{3}(250,85)(250,17)
\put(247,15){$\bullet$}
\put(247,5){$b$}
\put(250,20){\vector(1,0){30}}
\put(264,24){$v$}

%\qbezier(185,80)(175,45)(185,5)
%\dottedline{2}(88,270)(257,290)

         \end{picture}

      \caption{Manifolds $M,B$ and the map $\pi$ in
      {\bf Counterexample~\ref{lcex425}}.}
         \label{fig412cex425}
         \index{pictures!Counterexample\ref{lcex425}}
      \end{figure}
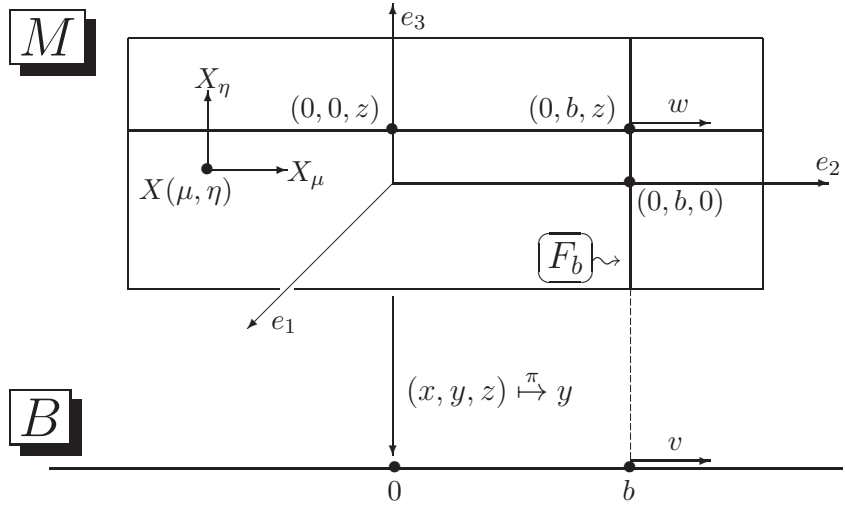

%(see Fig.~{\ref{fig412cex425}})

Each fiber is given by,
\begin{eqnarray*}
F_{b} & = & \pi^{-1}(b)=
\{(0,b,z): z\in\real\}, \hspace{0.1in}b\in B=\real
\end{eqnarray*}
which is a line passing through $(0,b,0)$, determined by the
intersection of $M$ with the plane $y=b$. Hence each fiber is
\underline{not compact} as a subset of \ra{3}, they all have  infinite
diameter, and therefore the fibers' diameters are
\underline{not uniformly bounded}.

Our goal next is to show that $\pi$ is not a rough isometry.

It suffices to verify that $\pi$ does not satisfy {\bf (RI.1)} for
particular pairs of elements in $M$, i.e.,

$\forall A\geq 1$, $\forall C>0$, 
$\exists \eta_{AC}\in\real\setminus\{0\}$, 
\begin{eqnarray}
\protect\label{leqn4412}
& & d_{B}\left(\pi\circ X(\mu,\eta),\pi\circ X(\mu,0)\right)<
\frac{1}{A}\cdot d_{M}\left( X(\mu,\eta),X(\mu,0)\right)-C
\Leftrightarrow\nonumber\\
& \Leftrightarrow & |\mu-\mu|<
\frac{1}{A}\cdot d_{M}\left((0,\mu,\eta),(0,\mu,0)\right)-C
\Leftrightarrow\nonumber\\
& \Leftrightarrow & 0<\frac{1}{A}\cdot|\eta|-C,
\hspace{0.1in}\forall \eta\geq\eta_{AC} 
%\hspace{0.5in}\ce{25}
\end{eqnarray}

Let $A\geq 1, C>0$ be arbitrary, and define the real positive 
number $\eta_{AC}:=AC+1>0$.
 
We see that,
\begin{eqnarray}
\protect\label{leqn4413}
\eta_{AC}=AC+1>AC \Leftrightarrow \frac{1}{A}\eta_{AC}>C
\Leftrightarrow \frac{1}{A}\eta_{AC}-C>0
%\hspace{0.5in}\ce{26}
\end{eqnarray}
and since, for all $\eta\geq\eta_{AC}$,
\[
\frac{1}{A}\eta-C\geq\frac{1}{A}\eta_{AC}-C
\stackrel{(\ref{leqn4413})}{>}0
\]
inequality (\ref{leqn4412}) is verified.

Therefore, (\ref{leqn4412}) 
holds and {\bf (RI.1)} fails for $\pi$, which shows that $\pi$ is
\underline{not} a rough isometry.

This describes the Counterexample.

\vspace{0.2in}

\begin{center}
{\bf Acknowledgments}
\end{center}

%I wish to express my gratitude to everyone who contributed 
%to making this article a reality.

I wish to register my sincere gratitude and thanks to
Professor Edgar Feldman, my doctoral thesis advisor, and to 
Professor Christina Sormani for assistance with exposition.

%\vspace{5in}

%\begin{flushright}
%{\it In memory of,}\hspace{0.9in}.\\
%{\it Mrs. Joana Reis-Abreu},\\
%\end{flushright}

  \nocite{*}

  \bibliography{maxrankarxiv}

%\vspace{1in}

\begin{flushright}
%C. Abreu-Suzuki\\
Mathematics, CIS Department\\
SUNY College at Old Westbury\\
Old Westbury, NY 11568\\
email: casuzuki@earthlink.net\\
\end{flushright}

\end{document}